  \documentclass{amsart}                 
 
 \usepackage{graphicx}
  \usepackage[format=plain,font=small]{caption}
 
\usepackage{amssymb,amsmath,amsfonts}

\def\d{\delta}

\def\R{{\mathbb R}}
\def\C{{\mathbb C}}

\def\Z{{\mathbb Z}}
\def\cM{{ \mathcal M}} 
\def\cL{{ \mathcal L}}

\def\cF{{ \mathcal F}}

\def\cR{{ \mathcal R}}
\def\cL{\mathcal{L}}

\def\o{{\omega}}

\newcommand\ove[1]{\mathop{\overline{#1}}}

 \def\snc{{\rm sinc}}
\newcommand\what[1]{\widehat{#1}}
\newcommand\essinf{\mathop{\rm ess\,inf}}

\newcommand\esssup{\mathop{\rm ess\,sup}}

\newcommand\leng[2]{{\rm len}_{\phantom{}_{#2}}({#1})}
\newcommand\rango[1]{{\rm rank\, }{{#1}}}

 \newcommand\cp{\mathop{{\text{{\,\large$\times\, $}}}}}
 \newcommand\pr{\mathop{{\text{\Huge{$\times$}}}}}
  \newcommand\crossp[2]{\mathop{\pr_{#1}^{#2}}}
\newcommand\crosspd[3]{\mathop{\pr_{#1}^{#2}}_{#3}}
\newcommand\crosspw[5]{\mathop{\pr_{#1}^{#2}}{#3}\,\big\langle\, {#4}\leftarrow{#5}\,  \big\rangle}

\def\Ld{{L^2(\R)}}  
\def\Lu{{L^1(\R)}}  
 \def\Lp#1#2{{L_{#1}^{#2}}}     
    
 \def\LdZ {{L_h^2(\R;\ell^2(\Z))}}

\def\ld{{\ell^2(\Z)}}

 \def\ldCN {{{\ell}^2 (\Z;{\C}^N)}}  

\def \Tfi{{T^{}_{\Phi,\hskip 0.05em t_o}}}
\def\Tfia{{T^{*}_{\Phi,\hskip 0.05em t_o}}}
\def\Tfistar{{T^{}_{\Phi^*,\hskip 0.05em t_o}}}
\def\Tfistara{{T^{*}_{\Phi^*,\hskip 0.05em t_o}}}

\def\Efi{{E_{\Phi,\hskip 0.05em t_o}}}
\def\Efistar{{E_{\Phi^*,\hskip 0.05em t_o}}}

\def\Sfi{{S_{\Phi,\hskip 0.05em t_o}}}
\def\BL{{B_\omega}}
\def\Gfi{{G_{\Phi,\hskip 0.05em t_o}}}
\def\Gfimu{{G^{-1}_{\Phi,\hskip 0.05em t_o}}}
\def\Gfitilde{{\widetilde{G}_{\Phi,\hskip 0.05em t_o}}}
\def\Gfitildecal{{ \widetilde{\mathcal{G}}_{\Phi,\hskip 0.05em t_o}}}
\def\Jfi{{J^{}_{\Phi,\hskip 0.05em t_o}}}
\def\Jfidual{{J^{}_{\Phi^*,\hskip 0.05em t_o}}}
 \def\Jfistar{{J^{*}_{\Phi,\hskip 0.05em t_o}}}
 
\newcommand\Apen[1]{{ {[\,\Jfi\,]}_{#1}}}
 \newcommand\Jbfidual{\mathbb{J}^{}_{\Phi^*,t_o}}
 \newcommand\Jbfistar{\mathbb{J}^{*}_{\Phi,t_o}} 
   
  \newcommand\Jbfi{\mathbb{J}_{\Phi,t_o}}
 \newcommand\JJbstinv{\mathop{(\det\,\Jbfi\Jbfistar )^{-1}}}
 
 \def\Ofi#1{{\Omega^{#1}_{\Phi,t_o}}}

\def\S-k{{ {S\mskip 5mu \hat{ }}\!\!_{-k}}}

\def\F^{{{ \hat  F} }}
\def\FF^{{ { \hat  F}^*{ \hat  F} }}

 \newcommand\Du{{ E}}

\def\squareforqed{\hbox{\rlap{$\sqcap$}$\sqcup$}}
\def\qed{\ifmmode\squareforqed\else{\unskip\nobreak\hfil
\penalty50\hskip1em\null\nobreak\hfil\squareforqed
\parfillskip=0pt\finalhyphendemerits=0\endgraf}\fi} 
 \newtheorem{theorem}{Theorem}[section] 
\newtheorem{corollary}[theorem]{Corollary}
\newtheorem{lemma}[theorem]{Lemma}

\newtheorem{proposition}[theorem]{Proposition}

\theoremstyle{remark}

\numberwithin{equation}{section}

\newcommand\rmi{\hbox{\rm (i)}}
\newcommand\rmii{\hbox{\rm (ii)}}

\begin{document}
  
\title[Frames for band limited functions]
{Frames and Oversampling Formulas\\for Band Limited  Functions}

\subjclass[2000]{} 


 \author[Vincenza del Prete]
 {Vincenza Del Prete}

 \address{Dipartimento di Matematica\\
 Universit\`a di Genova, via Dodecaneso 35, 16146 Genova \\ Italia}
 
\begin{abstract}
In this article we obtain families of  frames  for the  space $\BL$ of   functions with band in $[-\o,\o]$ by using  the theory   of  
  shift-invariant spaces. Our results are based on the Gramian analysis  of  A. Ron and  Z. Shen and a variant, due to Bownik, of their characterization of families of functions  whose shifts form  frames or Riesz bases.  We give necessary and sufficient conditions  for the translates of a finite number of functions (generators) to be a frame or a Riesz basis for $\BL.$ We also give explicit formulas for the dual generators
 and we apply them to Hilbert transform sampling and derivative sampling.  Finally, we   provide numerical experiments which support the theory.
\end{abstract}

\maketitle
\setcounter{section}{0}
\section{Introduction} \label{s:Introduction}
 
 In many signal and image processing applications, images and signals are assumed to be band limited. A band limited signal is a function which belongs to the  space $\BL$  of functions in $\Ld$ whose Fourier transforms have support in $[-\o,\o]$. Functions belonging to this   space can be represented  by  the Whittaker-Kotelnikov-Shannon series, which    is  the expansion   in terms of the  orthonormal basis of translates of the {\it sinc} function. The  coefficients
 of the expansion   are the  samples of the function  at  a uniform grid on $\R$,  with  ``density" $\frac{\o}{\pi}$. This sampling density is usually   called  the  {\it Nyquist} density.   \par\noindent The theory has been extended in many directions by several authors. In one of these extensions  the sinc orthonormal basis has been replaced   by  Riesz bases    (see for example    \cite{Hi} by J. R. Higgins) or frames. Frames   generally are overcomplete  and their expansion coefficients  are not unique. Their  redundancy is useful in applications  because  the reconstruction is more stable with respect to errors in the calculation of the coefficients 
 and it allows   the recovery of missing samples \cite{F}. In signal analysis  frames can be viewed as an ``oversampling" with respect to the Nyquist density.\par\noindent
 The second extension  consist in using more than one function  to generate the space $\BL$. In this case  the Riesz basis or the frame are formed by the translates of a  finite   family of  functions and  the expansion formula is called a multi-channel sampling formula   \cite{P}  \cite{Hi1}.\par 
Finally we mention that there is a huge literature on the problem of reconstruction of signals from non-uniform samples. Since this paper is only concerned with uniform sampling,  we refer  the reader to the recent article of A. Aldroubi and  K. Gr$\bf \rm{   \ddot{o}}$chenig 
 \cite{AG}   and the references given there for an extensive review of the problem of non-uniform sampling.
    \par
In this paper we construct multi-channel  uniform sampling formulas for band-limited functions using the theory of frames for shift-invariant spaces.
 A $t_o$-shift-invariant  space is a  subspace of $\Ld$ that is invariant   under all translations  $\tau_{k t_o}, k \in \Z$, by integer multiples of a positive number $t_o$.  We recall that   $\BL$ is  $t_o$-shift-invariant  
for any $t_o.$   A set  $\Phi $ in a $t_o$-shift-invariant space $S$     is called a set of generators  if $S$   is the closure of the space generated by the  family $\Efi=\{ \tau_{kt_o}\ \varphi,\,  \varphi\in \Phi,\,  k\in \Z \}$. 
 The space $S$  is said to be finitely generated if it has a finite set of generators. Finitely generated shift-invariant spaces can have different sets of generators;
the smallest number of generators is called the length of the space. \par\noindent
   The   structure of finitely generated shift-invariant spaces was investigated  by
C. de Boor, R. DeVore and A. Ron with the use of fiberization techniques based on the range function \cite{BDR}. 
These  authors  gave  conditions under which a finitely generated shift-invariant  space has a generating set satisfying 
  properties like  stability and  orthogonality. Successively, A. Ron and Z. Shen introduced  the  Gramian analysis  and extended the  results of \cite{BDR} to countable generated SI  spaces \cite{RS}. In  their  paper  they characterized sets of generators whose translates form    Bessel sequences, frames and Riesz bases. For finitely generated spaces these conditions are expressed in terms of the eigenvalues of the Gramian matrix.   In concrete cases it would be useful to have more explicit conditions  expressed in terms of the generators or their Fourier transforms.  In this paper, using a   result of  M. Bownik,  we obtain these more explicit conditions for  the space $\BL$ \cite{B}.  We also give explicit formulas for the Fourier transforms of the dual generators. \\
 In the last part of the paper  we use these results  to obtain     multichannel   oversampling formulas for band limited signals  and we apply them to the derivative sampling. Finally, we support the theory  with   some numerical experiments.   \par\noindent 
 The paper is organized as follows. In Section 2  we   collect some  results on frames for shift-invariant spaces and     find the range of $\BL$ as  a $t_o$-shift-invariant space.\\ 
In  Section 3 we consider a family $\Phi=\{\varphi_1,\varphi_2 , \ldots   \varphi_N\}$ of  $N$  generators, where $N$ is the length of ${\BL}$  as  a  $t_o$-shift-invariant space. 
 We give  a necessary and sufficient condition  for $\Efi$ to be a frame or a Riesz basis for $\BL$. The condition  is expressed in terms of the   pre-Gramian.   To prove this condition we use a result of Bownik which characterizes the system of translates $\Efi$ as being a frame or a  Riesz family  in terms of   ``fibers"   \cite{B}.  \par\noindent
 In  Section 4 we   find  the family $\Phi^*$  of dual generators. Here we use   the fact that,  in the fibered representation of $\BL$,  the frame   operator  is unitarily equivalent to the operator of multiplication by the dual Gramian matrix $\Gfitilde(x)=\Jfi(x)\Jfistar(x)$ acting on the fiber over $x.$
   The problem of finding the dual generators is reduced to solving   the  matricial equation $ \Jfi(x)=\Jfi(x)\Jfistar(x)\, \Jfidual(x)  $ in   the unknown  $\Jfidual(x)$    in $[0,h]$. Where $\Jfi(x)$   is invertible,  the solution   is the inverse of $\Jfi(x)$, elsewhere it is given by  the Moore-Penrose inverse  of $\Jfi(x)$. We  give the expressions of the Fourier transforms of the dual generators as cross   products of  translates     of  the vector   $\what{\Phi}=(\what{\varphi_1},\what{\varphi_2}, \ldots \what{ \varphi_N}).$\par\noindent
In  Section 5 we apply  the results of the previous sections  to obtain two and three-channel sampling formulas for   functions of    $\BL$ and we apply them  to   derivative sampling.  We also provide some numerical experiments. 
 
 \section{Preliminaries}
\label{s: Preliminaries}

In this section  we   collect some  results on frames for shift-invariant spaces   to be used later. We begin by introducing some notation. 
The Fourier transform of a  function  $f$ in 
  $  \Lu $ is
 $$ \cF f(\xi)={\hat f}(\xi)= \frac{1} {\sqrt{2\pi}} \int f(t)e^{-it\xi}dt.$$   

 The convolution of two functions $f$ and $g$ is
$$ f\ast g (x)=\frac{1} {\sqrt{2\pi} }\int f(x-y)g(y)\,dy,$$
so that $\cF (f\ast g) =\cF f\,\cF g.$
Let $h$ be a positive real number; $    \Lp{h}{p}$ is the space of $h$-periodic functions on $\R$ such that 
\begin{equation}\label{normaper} \|f\|_{\Lp{h}{p}}=\Bigl(\frac{1}{h}\int_0^h|f(x)|^p dx\Bigr)^{1/p}<\infty.
  \end{equation} 
 With the symbol $ \ldCN $ we shall denote the  space of $\C^N$-valued sequences $c=(c(n))_{\Z}$   such that
$$\|c\|_{\ell_2} =\Bigl(\sum_{n\in \Z} |c(n)|^2\Bigr)^{1/2} < \infty.$$  \par\noindent
 Let $H$ be a subspace of $\Ld.$   Given a subset $\Phi=\{ \varphi_j ,j=1,\dots, N\} $ of $H$ and a positive number $t_o$  denote by $ \Efi$ the set    $$ \Efi =\{\tau_{nt_o}\varphi_j , \ n\in \Z, \ j=1,\dots ,N\}.$$
  Here   $\tau_a f(x)=f(x+a).$   The closure   of the space generated by $ \Efi$ will be denoted   by $\Sfi$. 
 The family $ \Efi $ is   a frame for     $H$   if the operator $\Tfi: \ldCN \rightarrow  H$   defined by
 \begin{equation}\nonumber
\Tfi c=\sum_{j=1}^N \sum_{n\in Z} c_{j}(n) \tau_{n t_o}\varphi_j
\end{equation}
is     continuous,  surjective and ${\rm{ran}}(\Tfi)$ is  closed.   The family  $ \Efi
$ is a   frame  for $H$  if  and  only if there exist two constants $0<A\le B$ such that \begin{equation}\label{defframe}
A\|f\|^2\leq \sum_{j=1}^{N}\sum_{n\in \Z}  |\langle f,  \tau_{n t_o}\varphi_j \rangle  |^2
\leq B \|f\|^2 \hskip 1truecm \forall f\in  H.
\end{equation}   The constants $A$ and $B$  are called frame bounds. If $A=B$  the frame is called {\it tight }   and if $A=B=1$ a {\it Parseval} frame.   Denote  by   $\Tfia: H \rightarrow \ldCN$ the    adjoint of $\Tfi,$ defined by
\begin{equation}\label{Tstella}
(\Tfia f)_j(n)= \langle f,  \tau_{n t_o}\varphi_j  \rangle  \hskip1truecm n\in \Z ,  j=1\dots ,N.
\end{equation}
The operator  $\Tfi\Tfia:H\rightarrow H$  is called {\it frame operator.} 
 The set $\Efi$ is a frame for $H$ if and only if  the frame operator is continuously invertible and   \begin{eqnarray}\nonumber
\Tfi\Tfia  f = \sum_{j=1}^{N} \sum_{n\in \Z}  \langle f,  \tau_{n t_o}\varphi_j \rangle \tau_{n t_o}\varphi_j    \hskip1 truecm  f\in H.
\end{eqnarray} Observe that    
(\ref{defframe}) can be written  
$  A I\leq \ \Tfi \Tfia \  \leq BI, $
where $I$ is the identity operator on  $H.$ Denote  by $ \Phi^*$ the family  $\Phi^*=\{\varphi_j^*, j=1,\dots ,N\},$   where
\begin{equation}\label{fiduale}
 \varphi_j^*=(\Tfi\Tfia)^{-1}\varphi_j  \hskip1truecm     1\le j\le N.
 \end{equation}
If $\Efi$ is a frame for $H$ then    $\Efistar $ is also a frame
(the {\it dual frame}), and  
$\Tfi\Tfistara=\Tfistar\Tfia =~I$.  Explicitly
\begin{equation}\label{expansion}
f= \sum_{j=1}^{N} \sum_{n\in \Z}  \langle f,  \tau_{n t_o}\varphi_j^* \rangle \tau_{n t_o}\varphi_j= \sum_{j=1}^{N} \sum_{n\in \Z}  \langle f,  \tau_{n t_o}\varphi_j \rangle \tau_{n t_o}\varphi_j^* \hskip1truecm \forall f\in H.
 \end{equation}  
 The elements of  $\Phi$ are called {\it generators} and the elements of $\Phi^*$ {\it dual generators.} 
 If the family    $ \Efi
$ is a frame for $H$ and the operator $\Tfi$ is injective, then $ \Efi
 $ is called a  {\it  Riesz basis}. \\
 In what follows    $t_o$ is a positive parameter.  To simplify notation, throughout the paper we shall set  $$h=\frac{2\pi}{t_o}.$$
A subspace $S$ of $\Ld$   is $t_o$-shift-invariant if it is invariant under all  translations by  a multiple of $t_o.$  
The following bracket product  plays an important role in  Ron and Shen's analysis of shift-invariant spaces. For   $f$ and $g \in\Ld,$  define 
\begin{equation}\label{braket}
[f,g]=h\sum_{j\in \Z} f(\cdot+jh) \overline{g}(\cdot+jh).
\end{equation}
Note that $[f,g]$ is in $\Lp{h}{1}$ and $\|[f,f]\|_{\Lp{h}{1}}=\|f\|_2^2.$ 	 The  Fourier
coefficients of $ [\hat{f},\hat{g}]$ are given by
\begin{equation}\label{braketcoeff}
[\hat{f},\hat{g}]{\,}^{\widehat {\,}} (\ell) = \langle f, \tau_{\ell t_o}  g\rangle\qquad \ell\in \Z. \end{equation}
Indeed  $$
[\hat{f},\hat{g}]{\,}^{\widehat {\,}} (\ell) =\int_0^{h}
\sum_{j }\tau_{jh}(\hat{f}\, \overline{\hat{g}})(x) e^{-2\pi i  \ell \frac{x}{h}}dx 
 =\int \hat{f}(x)\, \overline{\hat{g} }(x)\  e^{-  i  \ell t_o x}dx$$
  If $S$ is a $t_o$-shift-invariant    space and  there exists     a finite family  $\Phi$ such that   $S=S_{\Phi,t_o},$ then we say that $S$ is finitely generated. Riesz bases for   finitely generated 
 shift-invariant spaces have been studied by various authors. In \cite{BDR} the authors  give a characterization of such bases.     A   characterization of frames and tight frames    also for countable sets $\Phi$ has been given by Ron and Shen in [RS1].   The principal notions of their theory are the {\it    pre-Gramian}, the {\it Gramian} and the {\it dual Gramian} matrices.
 \\
The pre-{\it Gramian} $\Jfi$ is the $h$-periodic function mapping $\R$ to the space of $\infty\times N$-matrices defined on $[0,h]$  by
\begin{equation}\label{pregramian}
\bigl(\Jfi \bigr)_{j\ell}(x)=\sqrt{h}\  \what{\varphi_\ell}(x+jh),\hskip1truecm j\in \Z, \,  \ell=1,\dots,N.
\end{equation}
The pre-Gramian $\Jfi$ should not be confused with the  matrix-valued function  whose entries  are $\sqrt{h}\, \what{\varphi_\ell}(x+jh),$ for all $x\in\R,$ which is not periodic. The spectrum of the space $\Sfi$ is defined as
\begin{equation}\label{spettro}  \sigma (\Sfi)=\{x\in \R: \Jfi(x)\not=0\} \end{equation}
or, equivalently,   as  the support of $\sum_{j=1}^N [\what{\varphi_j},\what{\varphi_j}].$
Of course, since  the functions $\what{\varphi_j}$   are defined only up to a null-set, the  support is intended  in the sense of distributions, i.e. as  the complement of the largest  
open set on which  the function $\sum_{j=1}^N [\what{\varphi_j},\what{\varphi_j}]$ vanishes as distribution. It was proved in \cite{BDR} that the spectrum of  a finitely generated space   depends only  on the space itself  and not on the particular  selection of its  generators.\\
  Denote by $\Jfistar$ the adjoint of $\Jfi.$ The {\it Gramian} matrix $\Gfi= \Jfistar \Jfi$ is the $N\times N$ matrix whose elements are the $h$-periodic functions  
  \begin{equation}\label{gramian}
\bigl(\Gfi \bigr)_{j\ell}=[\what{\varphi_{\ell}}  ,\what{\varphi_j}].
\end{equation}
  The    {\it dual Gramian} $\Gfitilde=\Jfi \Jfistar$ is the    infinite matrix whose elements are
\begin{equation}\label{dualgramian}
\bigl(\Gfitilde  \bigr)_{j\ell}=h \sum_{n=1}^N
\tau_{jh}\what{\varphi_{n}} \  \tau_{{\ell}h}\overline{\what{\varphi_{n}}} ,\hskip1truecm j, \ell\in \Z.
\end{equation}
The importance of these two  matrices lies in the fact that  the  Gramian matrix  represents the operator  $ \Tfia \Tfi$ 
 and the dual Gramian  represents  the operator  $ \Tfi   \Tfia$ and  many    properties of these operators  can be studied  by looking at them.  Indeed, by the  theory of Ron and Shen  \cite{RS}, after  conjugating  with  an isometry, the operator $\Tfi\Tfia$ can be decomposed into a measurable field of operators, acting on $\ld,$ which are represented by the  dual Gramian matrix $\Gfitilde$ in the canonical basis (see formulas (\ref{LFTTsF}) and (\ref{calMMfi}) below).
 Similarly,  after conjugation with   a Fourier transform, the  operator $\Tfia\Tfi$ is represented by the Gramian $\Gfi.$\\
Denote by $\LdZ$ the Hilbert space of $\ld$-valued $h$-periodic  functions  on $\R$ such that 
\begin{equation}\label{elleduevalued}
\|f\|_{\LdZ}=\Big(\frac{1}{h}\int_0^h\|f(x)\|_{\ell^2}\ dx\Big)^{\frac{1}{2}}<\infty.
\end{equation}
For every $f$ in  $\Ld$ we denote by $\cL_h f$ the $\ld$-valued function  defined on $[0,h]$ by 
\begin{equation}\label{opelleh}
 \cL_h f(x)  =\sqrt{h}\sum_{l\in \Z}f(x+\ell h)\delta_{\ell} \hskip 1truecm  x\in[0,h] 
 \end{equation}
and extended to $\R$ as a periodic function of period $h.$ Here $\{\delta_{\ell}: \ell\in \Z\}$ is  the canonical basis of $\ld.$ The map $f\mapsto \cL_hf$ is an isometry of $\Ld$ onto  $\LdZ.$ 
Observe  that   the vectors  $\cL_h \widehat{\varphi_j},\   j=1,\dots,N$
are the columns of the pre-Gramian $\Jfi$, i.e.
\begin{equation}\label{colonnepreg}
\Jfi=\big(\cL_h \widehat{\varphi_1},\ldots,\cL_h \widehat{\varphi_N}\big).
\end{equation}
  The   map $\cL_h$     links $t_o$-shift-invariant subspaces  of $\Ld$ with $h$-doubly-invariant subspaces of $\LdZ.$
  We recall that a  subspace of $\LdZ$ is  $h$-doubly-invariant if it is invariant  under pointwise multiplication by $e^{2\pi ik \frac{x}{h}},k\in \Z.$  Obviously a subspace $S$ of $\Ld$ is $t_o$-shift-invariant if and only if  the space  \begin{equation}\label{emmediesse}
 \cL_{h}(\hat{S})=\{\cL_h\hat{f}, f\in S \}
\end{equation}
  is $h$-doubly-invariant. T.P. Srinivanasan gave a  characterization of doubly-invariant spaces in terms of {\it range functions}  {(see \cite{H},\cite{S})}.  As remarked by  de Boor,  DeVore and  Ron   \cite{BDR}  a similar characterization  of shift-invariant spaces follows from it (see Proposition  \ref{CondRange} below).
 In our context  a  range function  is a $h$-periodic map $\cR$ from $\R$ to the closed subspaces  of $\ld.$ The function $\cR$ is measurable if the map  $\mathcal P$  which maps  a point $x\in \R$ to the orthogonal projection $\mathcal P(x)$  onto $\cR(x) $ is weakly measurable in the operator sense, i.e.     the function
$$x\mapsto (P(x)\varphi,\psi)_{\ld}$$
is measurable for all $\varphi$ and $\psi\in \ld.$   Range functions which are equal almost everywhere are identified. \begin{proposition} \label{CondRange} 
A  closed subspace $S$ of $\Ld$ is $t_o$-shift-invariant if and only if   
\begin{equation} 
S=\{f\in \Ld : \cL_h \hat{f}(x)\in \cR(x)\  {\rm for}\  a. e.\, x\in\R\} \end{equation}
for some measurable $h$-periodic range function $\cR.$ \end{proposition}
Obviously $t_o$-shift-invariant subspaces with the same range function coincide. This observation shall be used in the proof of  Theorems \ref{maingen} and \ref{maingendue}.\\
 In Theorem \ref{noteTeoX}  below we   compute the range function of the   space of   band-limited  functions 
 $$\BL=\{ f\in L^2(\R): supp( \hat{f})\subset [-\o,\o] \}.$$ 
 In general it is not a simple matter to compute the range function of a space.  Fortunately, if the space $S$ is finitely generated, the range function can be written in terms of the generators; indeed in \cite{BDR} it was shown that if $S=\Sfi$ and $\Phi=\{\varphi_j,\quad j=1,\dots N\}$ then
\begin{equation}\label{rangephi}
\cR_h(\Sfi)(x)= {\rm span}\{ \cL_h \what{\varphi_j}(x): j=1,\dots, N\} 
\end{equation}
is the space generated by  the columns of the pre-Gramian matrix $\Jfi$. This result holds also for  countable sets of generators.  
Of course  a shift-invariant space  $S$ can have more than one family of generators; the smallest number of generators is called the length of the space  
\begin{equation}\label{len}
\leng{S}{t_o}={\rm min}\{\sharp \Phi, S=\Sfi \}.\end{equation} 
 In \cite{ BDR}  it has been proved that  if $S$ is finitely generated then\begin{equation}\label{lenEds}\leng{S}{t_o}=\esssup\big\{ {\rm dim}\cR_h(\Sfi)(x):\quad  x\in [0,h]\big\}. 
\end{equation}
To state the next  theorem we   need some notation. We denote by $[a]$ the greatest integer less than $a.$ 
 Recall that $h= 2\pi/{t_o}$  and  set $\ell=\left [ \frac{\o}{h} \right ]+1.$ \par
If  $\frac{\o}{\ell}\le h< \frac{\o}{\ell-\frac{1}{2}}$ then
 $0\le -\o+\ell h<\o-(\ell-1)h<h$. We denote by
$I_{-},I,I_{+}$ the intervals defined by
  \begin{equation}\label{I}
I_{-}=(0,-\o+\ell h),\qquad  I= (-\o+\ell h,\o-(\ell-1)h), \qquad  I_+=(\o-(\ell-1)h,h) 
\end{equation} 
 Similarly  if   $\frac{\o}{\ell-\frac{1}{2}}\le h< \frac{\o}{\ell-1} $ then
 $0<\o- (\ell-1)h  \le -\o+\ell h<h;$ in this case  we denote by  $K_{-},K,K_{+}$ the intervals defined by
\begin{equation}\label{K}
K_{-}=(0,\o-(\ell-1) h),\qquad  K=  (\o-(\ell-1) h,-\o+\ell h),\qquad  K_+=(-\o+\ell h,h).
\end{equation}
Recall that  $\{\delta_{\ell}: \ell\in \Z\}$ is  the canonical basis of $\ld.$
 \begin{theorem}\label{noteTeoX} Let $t_o$ be a positive parameter,   $h= 2\pi/{t_o}   $  and  set $\ell=\left [ \frac{\o}{h} \right ]+1.$   Then 
 \begin{itemize} 
  \item[\rmi] if $\frac{\o}{\ell}\le h<\frac{\o}{\ell-\frac{1}{2}}$  the range function of the space $\BL$ is  
 \begin{equation}\label{RangeterzoX}
 \cR_h({\BL })(x)   =
\begin{cases}
{\rm span }\{\delta_j:\    -(\ell-1)\le j \le \ell-1 \},&      x\in I_{-}\\
{\rm span }\{\delta_j:\  -\ell \le j\le \ell-1  \},&     x\in I\\
{\rm span }\{\delta_j:\  -\ell\le j\le \ell-2  \},&      x\in I_{+}.\\
\end{cases}
\end{equation}
  Note that if   $h=\frac{\o}{\ell}$ then  the     intervals $I_{-}$ and $I_{+}$ are empty. 
\item[\rmii] If $\frac{\o}{\ell-\frac{1}{2}}\le h< \frac{\o}{\ell-1} $ the range function of the space $\BL$ is  
\begin{equation}\label{Rangeterzo1X}
 \cR_h({\BL })(x)   =
 \begin{cases}
{\rm span }\{\delta_j:\   -(\ell-1)\le j \le \ell-1  \},&    x\in K_{-}\\
{\rm span }\{\delta_j:\  -(\ell-1)\le j\le \ell-2  \},&   x\in K\\
{\rm span }\{\delta_j:\  -\ell\le j\le \ell-2  \}, &    x\in K_{+}.\\
\end{cases}
 \end{equation}
Note that if $h=\frac{\o}{\ell-\frac{1}{2}}$ then the   interval $K$ is empty. 
\end{itemize} 
 \end{theorem}  
  \begin{proof} 
For the sake of simplicity we prove the theorem only for $\ell=2,$ i.e. 
$\o/2\le h<\o.$ The proof in the other cases is similar.
Denote by ${\mathcal  M}$ the space of all functions $g$ in $\LdZ$ such that $g(x)\in {\mathcal R}(\BL)(x)$ for a.e. $x\in[0,h].$
  Let ${\mathbf Q}$ denote the orthogonal projection of $\LdZ$ onto $ {\mathcal  M}.$
 By \cite[Lemma  p. 58]{H} for $g\in \LdZ$  we have 
 \begin{equation}\label{rango:1}
{\bf Q}g(x)=Q(x)g(x)\hskip 1truecm    a.e. \ x
\end{equation} 
where $Q(x)$ is the orthogonal projection of $\ld$ onto ${\mathcal R}(\BL)(x).$
Thus, to determine  ${\mathcal R}(\BL)(x)$  we only need to describe the projection $Q(x)$ for a.e. $x.$ 
Denote by  $\Lambda_{\o}$   the  space of   functions in $\Ld$  with support in  $[-\o,\o]$ and   by ${\mathbf P}:\Ld\mapsto \Lambda_{\o}$
 the orthogonal projection onto  $\Lambda_{\o},$ that is ${\mathbf P}f=\chi_{[-\o,\o]}f.$ By Proposition \ref{CondRange}  $\cL_h(\Lambda_{\o})=\cL_h(\widehat{\BL})={\mathcal M}.$
Thus we have \begin{equation}\label{rango:2}{\mathbf{Q}}=\cL_h {\mathbf P} \cL_h^{-1} 
\end{equation}
because $\cL_h $ is an isometry.
Let $\psi=\sum_{n\in\Z}\psi_n \delta_n \in \LdZ,$ then 
 $ \cL_h^{-1} \psi(x)=\psi_n(x)$  for $x\in[nh,(n+1)h],$ $n\in \Z;$  
i.e.  \begin{equation}\nonumber
\cL_h^{-1} \psi = \sum_n \chi_{[nh,(n+1)h]} {\psi}_n. \end{equation}
Hence
\begin{equation}\label{PLhinvterzo}
{\mathbf P}\cL_h^{-1} \psi  =\sum_n \chi_{[-\o,\o]}\chi_{[nh,(n+1)h]}  {\psi}_n.
\end{equation}
Since $\frac{\o}{2}\le h < \o$
\begin{equation}\nonumber 
\chi_{[-\o,\o]}\chi_{[nh,(n+1)h]}\,=\,\begin{cases}
\chi_{[0,h]} & {\rm if}\quad   n=0 \\
\chi_{[-h,0]}   & {\rm if} \quad n=-1\\
\chi_{[h,\o]}   & {\rm if} \quad n=1\\
\chi_{[-\o,-h]}  & {\rm if} \quad n=-2 \\
0 & {\rm otherwise}.
\end{cases}
\end{equation}
It follows that
 \begin{equation}\label{PLhinvancora}
{\mathbf P}\cL_h^{-1} \psi=\chi_{[-\o,-h]} \, {\psi}_{-2}+\chi_{[-h,0]} \, {\psi}_{-1}+
\chi_{[0,h]}{\psi}_0 +
\chi_{[h,\o]} \, {\psi}_{1}.
 \end{equation}
Now we     find $\cL_h {\mathbf P}\cL_h^{-1}\psi(x). $\\
First suppose  that    $ h \ge\frac{2}{3}\o;$ then
  $0<\o-h\le2h-\o< h.$   Hence from (\ref{PLhinvancora}) it follows that 
  for each $\psi\in \ld$
 \begin{equation}\nonumber
\cL_h{\mathbf P}\cL_h^{-1} \psi(x)=
\begin{cases}
  {\psi}_{-1}(x)\,  {\d}_{-1}+\psi_0(x)\,\d_0+\psi_{1}(x)\, \d_{1}   & \quad x\in \,  (0,\o-h) \\
 {\psi}_{-1}(x)\, \d_{-1}+{\psi}_0(x)\,  \d_0     & \quad        x\in \, (\o-h, 2h-\o)\\
{\psi}_{-2}(x)\, \d_{-2}+{\psi}_{-1}(x)\, \d_{-1}+ \psi_0(x)\,\d_0   & \quad x\in   \,    (2h- \o,h).\\
\end{cases}
\end{equation}
Suppose now  $ h <\frac{2}{3}\o;$ 
in this case $2h-\o<\o-h<h.$   Hence from (\ref{PLhinvancora})
we have
  $\psi\in \ld$
   \begin{equation}\nonumber
\cL_h {\mathbf P}\cL_h^{-1} \psi(x)=
\begin{cases}
  {\psi}_{-1}(x)\,  {\d}_{-1}+\psi_0(x)\,\d_0+\psi_{1}(x)\, \d_{1}   &   x\in \,    (0,2h-\o) \\ \displaystyle{\sum_{j=-2}^{1}
  {\psi}_{j}(x) \,\d_{j}  }& x\in \,     (2h-\o,\o-h)\\
  {\psi}_{-2}(x)\,\d_{-2}+{\psi}_{-1}(x)\,\d_{-1}+ \psi_0(x)\,\d_{0}    &  x\in \,     (\o-h,h).\\
\end{cases}
\end{equation}
Note that if $h=\o/2$ the intervals $(0,2h-\o)$ and $(\o-h,h)$ are empty. This determines completely the projection $Q(x)$ by (\ref{rango:1}) and (\ref{rango:2}).
  \end{proof}
  The following corollary is a straightforward consequence of   Theorem \ref{noteTeoX}.
  \begin{corollary} \label{lenBL} Let   $\ell=\left [ \frac{\o}{h} \right ]+1.$ Then the length of  $\BL$, as  $t_o$-shift-invariant space, is
\begin{equation}
\leng{\BL}{t_o}=\begin{cases}
2\ell   &   {\rm {if}}\ \frac{\o}{\ell}\le h<\frac{\o}{\ell-\frac{1}{2}}, \\
2\ell-1   & {\rm {if}} \   \frac{\o}{\ell-\frac{1}{2}}\le h< \frac{\o}{\ell-1}.\end{cases}
\end{equation} 
  \end{corollary}
  \begin{proof}
 The thesis follows immediately from (\ref{lenEds})   and Theorem \ref{noteTeoX}.     \end{proof}
 \section{Frames for the space $\BL $}
\label{s: Frames}
Let $\Phi$ be a finite family of generators for $\BL$. In this section   we find conditions under which $\Efi$ is a frame for $\BL$. First  we prove a representation formula for the frame operator analogous to the formula proved by Heil and Walnut for Gabor frames \cite[Theorem 4.2.1]{HW}.  
From this formula we deduce a simple necessary condition (see Proposition \ref{Prop:nec}).
 The formula will also be useful in Section \ref{s: Dual generators}  to find the dual generators.\par   
Let $\Phi=\{\varphi_j , 1\le j\le N\} $ be    a family  of functions in  $\BL,$   $t_o$ a positive parameter and $h={2\pi}/{t_o}.$ Denote by $\Ofi{k}$ the  function
   \begin{equation}\label{Omegacappa}
\Ofi{k}=h\sum_{j=1}^N \ \what{\varphi_j}{\tau_{kh}\overline{\what{\varphi_j}}}\hskip 1.5 truecm k\in \Z.
\end{equation}
 \begin{theorem}\label{HW}     
 If $\sum_k \| \Ofi{k}\|_{\infty}< \infty$ then the operator $\Tfi$ from $\ldCN$ to $\BL$ is bounded  and  \begin{equation}\label{HeilWalnut}
{ \mathcal F} \Tfi \Tfia { \mathcal F}^{-1} =\sum_k \Ofi{k} \tau_{kh}
\end{equation}
on $\Ld, $ where the series converges in the operator norm.\end{theorem}
 \begin{proof} 
First we show that if $f\in \Ld$  then $[\what{f},\what{\varphi_j}]$ is in 
$\Lp{h}{2} $ for  $1\le j \le N.$  Indeed 
 $$\begin{aligned}
   \int_0^h \big|[\what{f},\what{\varphi_j}] \big|^2 dx& =
 h \int_0^h\Big( \sum_{\ell}
    \overline{\what{f}\,}(x+\ell h) \,  \what{ \varphi_j}(x+\ell h)\Big)\,
  [\what{f},\what{\varphi_j}](x)  dx\\
 &=h \sum_{\ell}\int_{\ell h}^{(\ell+1)h} \overline{ \what{f}\,}(z) \,\what{\varphi_j}(z) \, [\what{f},\what{\varphi_j}](z) \, dz \\
&=h^2 \int \sum_k\what{\varphi_j}(z) \,\overline{\what{\varphi}}(z+kh) \,\what{f}(z+kh)\overline{\what{f}\,}(z)\,dz. 
\end{aligned} $$
Note that we may   exchange   the sum and the integral  because  $ \overline{ \widehat {f}\, }  \what{ \varphi_j}$  has compact support and the sum is finite. 
By summing over $j$ and exchanging the sums we obtain
$$
\sum_{j=1}^N  \frac{1}{h}  \int_0^h \big|[\hat{f},\what{\varphi_j}] \big|^2 dx=
  h \int \sum_k \sum_{j=1}^N \what{\varphi_j}(z) \overline{\what{\varphi}}(z+kh) \what{f}(z+kh)\overline{\what{f}\, }(z)\,dz. 
$$
Therefore, by (\ref{Omegacappa})
\begin{equation}\label{intermedia1}
\sum_{j=1}^N  \frac{1}{h}  \int_0^h \big|[\what{f},\what{\varphi_j}] \big|^2 dx = 
 \langle  \sum_k \Ofi{k} \tau_{kh} \what{f},\what{f}\rangle.
\end{equation} 
By  Schwarz's inequality   the right hand side  is less then  $  \|\hat{f}\|^2_2 \sum_k\| \Ofi{k}\|_{\infty} .$ Hence     $  [\what{f},\hat{\varphi_j}]\in  \Lp{h}{2},$ for $1 \le j \le N.$ By (\ref{braketcoeff}) the Fourier coefficients of 
$[\what{f},\what{\varphi_j}]$ are  $\langle f,\tau_{kh}\varphi_j \rangle,\  k\in \Z.$
Hence, by  (\ref{intermedia1}),  Plancherel's formula and   (\ref{Tstella}) 
\begin{align*} 
 \langle  \sum_k \Ofi{k} \tau_{kh} \what{f},\what{f}\rangle &=\sum_{j=1}^N \sum_k
 |\langle f,\tau_{kh}\varphi_j \rangle |^2\\ &=\|\Tfia f\|^2\\ &= 
 \langle \Tfi\Tfia f,f\rangle.
\end{align*} 
This  proves  (\ref{HeilWalnut}). Hence $\Tfia$ is bounded  from $\BL$ to \,$\ldCN$  and $\Tfi$  is bounded from $\ldCN$ to $\BL$.
 \end{proof} Now by using formula (\ref{HeilWalnut})
we   show that the operator  $\Tfi \Tfia$  is unitarily  equivalent  to the operator of multiplication by  the matrix  $\Gfitilde$ acting on the space $\LdZ.$  Our proof, given for the sake of completeness,  is a simple alternative derivation of a result  of Ron and Shen for general shift-invariant spaces \cite{RS}. \par\noindent 
 Let  $f$ be a function in $\Ld.$ By  (\ref{HeilWalnut})  for each $j\in \Z$
$$\begin{aligned}   
{ \mathcal F} \Tfi \Tfia { \mathcal F}^{-1} f(x+jh)  &=  \sum_k  \Ofi{k}(x+jh) f(x+jh+kh)\\
 &=  \sum_{\ell} \Ofi{\ell-j}(x+jh) f(x+\ell h). 
\end{aligned}
$$
Now we observe that 
$$\Ofi{\ell-j}(x+jh)=\overline{\Ofi{j-\ell}}(x+\ell h) \hskip1truecm {\rm for\ a.e. \ } x\in\R\ \  j,\ell \in \Z. $$
Therefore by (\ref{dualgramian}) and  (\ref{Omegacappa}) we get
\begin{equation}\nonumber
{\mathcal F} \Tfi \Tfia { \mathcal F}^{-1} \,f(x+j h)=
\sum_{\ell}\bigl(\Gfitilde  \bigr)_{j\ell}(x)\   f(x+\ell h)\hskip 2 truecm j,\ell \in \Z.
\end{equation}
By   (\ref{opelleh}) this formula implies that
\begin{equation}\label{LFTTsF}
\cL_h  {\mathcal F} \Tfi \Tfia {\mathcal F}^{-1} f=
\Gfitilde   \cL_h f  \hskip 1truecm [0,h].
\end{equation}
 This shows that   the operator $\Tfi \Tfia$ is unitarily equivalent to the
  operator $ \Gfitildecal $  defined by
 \begin{equation}\label{calMMfi}
 {  \Gfitildecal} g(x)=   \Gfitilde(x) g(x) \end{equation} 
  for almost every $x\in[0,h].$
 An operator of this form is  said to be  decomposable into the measurable field 
 $x\mapsto   \Gfitildecal(x)$ of operators on $ l^2(\Z).$ We shall use this representation in Section 3 to find the dual generators of frames. 
\par
In the  following proposition we give a necessary condition for $\Efi$ to be a frame. The proof mimics closely an argument of Heil and Walnut for  Gabor frames \cite{HW}. 
\begin{proposition}\label{Prop:nec} Let $\Phi=\{\varphi_j , 1\le j\le N\} $ be a family of functions of  $\BL.$ If  $\Efi$ is a frame for $\BL$ then 
there exist   $0<\delta\le\gamma<\infty $ such that 
\begin{equation}\label{main:1} 
 \delta\le    \sum_{j=1}^N |\what{\varphi_j}|^2\le \gamma \hskip 1truecm {\rm a.e. \   in }\ [-\o,\o]. \end{equation}
\end{proposition}
\begin{proof} Observe that $  h \sum_{j=1}^N |\what{\varphi_j}|^2 =\Ofi{0}.$ We only prove that  $0<\essinf \Ofi{0}$ because the proof of  the other inequality  is analogous.  Suppose that  $\essinf \Ofi{0}=0$; \break then  for each $\epsilon$ there exists a set $E_{\epsilon}\subset(\o,\o),$  of positive measure, such that $\Ofi{0}(x)<\epsilon $ for a.e.
 $x\in E_{\epsilon}.$ We may suppose that there exists an   an interval  $I$ of measure $h$ such that 
$E_{\epsilon}\subset I.$ Let $f\in \BL$ be defined by $\hat{f}=\chi_{{ }_{E_{\epsilon}}};$   by the Parseval and  the Plancherel  formula 
$$\sum_k|\langle f,\tau_{kt_o}\varphi_j\rangle|^2=\sum_k  \Big|
\int_{I}\what{f}\ \overline{\what{\varphi_j}} e^{2\pi i k\frac{x}{h}} \, dx\Big|^2=h\int_I \chi_{{ }_{E_{\epsilon}}}\big| \what{\varphi_j} \big|^2\ dx\hskip 1truecm j=1,\dots,N.$$  
Therefore 
$$\sum_{j=1}^N\sum_k|\langle f,\tau_{kt_o}\varphi_j\rangle|^2=
 \int_I\chi_{{ }_{E_{\epsilon}}}(x) \Ofi{0}(x)\ dx < \epsilon h \|\chi_{{ }_{E_{\epsilon}}} \|^2= \epsilon h \|f\|^2.$$ This contraddicts the  fact that $\Efi$ is a frame. \end{proof}
  Next, we give a  necessary and sufficient conditions  for 
 $\Efi$  to be  a frame or a Riesz basis for $\BL$, when  
  $\Phi$ is a subset of $\BL$ of cardinality $\leng{\BL}{t_o}$   (see Theorems     \ref{maingen} and   \ref{maingendue} below).   
 Our characterization will be given in terms of the pre-Gramian. Strictly speaking, the pre-Gramian $J_{\Phi,to}$ is an infinite matrix. However we shall see that all but a finite number of the rows of $J_{\Phi,to}$ vanish. Hence we may identify it with a finite matrix. We shall need the following 
 \begin{lemma}\label{pretraslate} Let $h$ be a positive number,   
  $\ell=\left [ \frac{\o}{h} \right ]+1$  and   $g\in \BL$. Let $I_{-},I,I_{+}$ and $K_{-},K,K_{+}$
be  the intervals defined in the previous section (see  (\ref{I}) and (\ref{K})). \begin{itemize} 
  \item[\rmi] Suppose that   $\frac{\o}{\ell}\le h<\frac{\o}{\ell-\frac{1}{2}}  $. Then  $\tau_{jh}\,g(x)=0$   
 \begin{equation}\label{pretraslateuno}
 \begin{aligned} 
if \   x\in I_{-}   \quad and \ &  j \, \notin \{-(\ell-1)\le j \le \ell-1\}     \\
if \   x\in I_{\ }\,  \quad and \ &  j\, \notin \{ -\ell \le j\le \ell-1\}     \\
 if\    x\in I_{+} \quad and \ &  j\, \notin \{ -\ell\le j\le \ell-2 \}.\\
\end{aligned}
\end{equation}
 Note that if   $h=\frac{\o}{\ell}$ then  $I_{-}$ and $I_+$  are empty. 
\item[\rmii] Suppose that  $\frac{\o}{\ell-\frac{1}{2}}\le h< \frac{\o}{\ell-1} $. Then  $\tau_{jh}\,g(x)=0$    
\begin{equation}\label{pretraslate2} \begin{aligned}
if \   x\in K_{-}   \quad and\  &  j\, \notin \{  -(\ell-1)\le j\le \ell-1\}\\
if \   x\in K_{\ }\,  \quad and \ &  j\, \notin \{-(\ell-1)\le j\le \ell-2 \}\\
 if\    x\in K_{+} \quad and \ &  j \, \notin \{-\ell\le j\le \ell-2\}.\\
\end{aligned}\end{equation}
\end{itemize}
 Note that if   $h=\frac{\o}{\ell-\frac{1}{2}}$ then    $K$ is empty. 
\end{lemma}
We omit the proof which   is straightforward.
\\
\par We consider separately the two cases $ h< \frac{\o}{\ell-\frac{1}{2}}$
and $ \frac{\o}{\ell-\frac{1}{2}}\le h.$
 Assume first that $\frac{\o}{\ell}\le h< \frac{\o}{\ell-\frac{1}{2}}.$ Then   $\leng{\BL}{t_o}=2\ell$ by Corollary \ref{lenBL}. 
Let    $\Phi=\{\varphi_j:1\le j\le 2\ell\}$ be a subset of $\BL$  of cardinality $2\ell$. By Lemma \ref{pretraslate}  all the rows of the matrix $\Jfi $ vanish except possibly $\big(\tau_{j h}\what{\varphi}_1,\tau_{j h}\what{\varphi}_2,\dots,\tau_{j h}\what{\varphi}_{2\ell}\big)$, $-\ell\le j\le\ell-1$. 
Thus we identify the infinite matrices $\Jfi$,  $\Jfistar$ and $\Gfi$ with their $2\ell\times2\ell$ submatrices corresponding to these rows. The  entries of the Gramian matrix are
 \begin{equation}\nonumber
\bigl(\Gfi \bigr)_{j k}=[\what{\varphi_{k}}  ,\what{\varphi_j}]\hskip 1.5truecm  1\le j\le 2\ell  \qquad 1\le k \le 2 \ell .
\end{equation}
By Lemma \ref{pretraslate} the $i$-th column  of $\Jfi,$   $1\le i\le 2\ell$ is 
\begin{equation}\label{colonne:0}
\sqrt{h}\begin{bmatrix} 0\\
\tau_{-(\ell-1)h}\what{\varphi_i}\\
\vdots \\
\what{\varphi_i}\\
\\ \vdots \\
\tau_{(\ell-2)h}\what{\varphi_i}\\
\tau_{(\ell-1)h}\what{\varphi_i}\\
\end{bmatrix}{\rm{in}}\  I_{-}\hskip0.9truecm
\sqrt{h}\begin{bmatrix} \tau_{-\ell h}\what{\varphi_i}\\
\tau_{-(\ell-1)h}\what{\varphi_i}\\
\vdots \\
\what{\varphi_i}\\
\\ \vdots \\
\tau_{(\ell-2)h}\what{\varphi_i}\\
\tau_{(\ell-1)h}\what{\varphi_i}\\
\end{bmatrix}\ {\rm{in}}\  I\hskip0.9truecm
\sqrt{h}\begin{bmatrix} \tau_{-\ell h}\what{\varphi_i}\\
\tau_{-(\ell-1)h}\what{\varphi_i}\\
\vdots \\
\what{\varphi_i}\\
\\ \vdots \\
\tau_{(\ell-2)h}\what{\varphi_i}\\
0
\end{bmatrix}\  {\rm{in}}\  I_{+} \hskip.05em .
\end{equation}
 We   note that  
  \begin{equation}\label{ranKJeG}   
 \rango{ \Gfi}= \rango{\Jfi}.
\end{equation}
   We shall use the following result of Bownik \cite[Thm 2.3]{B}  which characterizes the system of translates $E_{\Phi,t_o}$ as being  a frame or a Riesz family (for the space it generates) in terms of the ``fibers" $\{\cL_h\hat\varphi(x): \varphi\in \Phi\}$.  
\begin{theorem}\label{t: Bownik}
Suppose $\Phi\subset L^2(\R^n)$ is countable and let $H$ be the subspace of $L^2(\R)$ generated by $E_{\Phi,t_o}$. Then
\begin{itemize}
\item[\rmi] $E_{\Phi,t_o}$ is a frame for $H$ with constants $A, B$  if and only if $\{\cL_h\hat\varphi(x): \varphi\in \Phi\}$ is a frame for $\cR_h(H)(x)$ with constants $A,
B$ for a.e. $x\in [0,h]$.
\item[\rmii] $E_{\Phi,t_o}$ is a Riesz basis for $H$ with constants $A, B$  if and only if $\{\cL_h\hat\varphi: \varphi\in \Phi\}$ is a Riesz basis for $\cR_h(H)(x)$ with constants $A,
B$ for a.e. $x\in [0,h]$.
\end{itemize}
\end{theorem}
To apply Bownik's theorem in our context we need a simple lemma of linear algebra. Let $J$  be a $n\times m$ matrix with complex entries, $n\le m$; we shall denote by $\|J \|$ the norm of $J$ as linear operator from $\C^m$ to $\C^n$ and by $[J]_n$ the sum of the squares of the absolute values of the minors of order $n$ of $J.$
\begin{lemma}\label{simple lemma} Let $v_1,\ldots,v_m$ be $m$ vectors in $\C^n$, $m\ge n$, and denote by $J$ the matrix $(v_1,\ldots,v_m)$ whose $j$-th column is the vector $v_j$. 
\begin{itemize}
\item[\rmi] If $[J]_n>0$ then $\{v_1,\ldots,v_m\}$ is a frame of $\C^n$ with frame constants $A\ge [J]_n\,\Vert J\Vert^{2(1-n)}$, $B=\Vert J\Vert^2$. Conversely, if $\{v_1,\ldots,v_m\}$ is a frame for $\C^n$ with constants $A$ and $B,$ then $[J]_n\ge A^n$ and   $\Vert J\Vert\le B^{1/2}$.
\item[\rmii] If $m=n$ and $\det J>0$ then $\{v_1,\ldots,v_m\}$ is a Riesz basis of $\C^n$ with constants $A=\det(J)^2\,\Vert J\Vert^{2(1-n)}$ and $B=\Vert J\Vert^2$. Conversely, if  $\{v_1,\ldots,v_m\}$ is a Riesz basis for $\C^n$ with constants $A$ and $B$ then $\det(J)\ge A^{n/2}$ and $\Vert J\Vert\le B^{1/2}$. \par
\end{itemize}
\end{lemma}
\begin{proof} 
Let $T:\C^m\to\C^n$ be the synthesis operator associated to $\{v_1,\ldots,v_m\}$, i.e.  $Tz=\sum_{i=1}^m z_j v_j$ for all $z\in\C^m$.  We observe that $TT^*$ and $JJ^*$ have the same eigenvalues $\lambda_1\le\cdots\le \lambda_n$, because the matrix $J$ represents the operator $T$ with respect to the canonical bases of $\C^m$ and $\C^n$. Since
$$
\lambda_1\,I \le TT^*\le \lambda_n \,I
$$
 $\{v_1,\ldots,v_m\}$ is a frame for $\C^n$  if and only if $\lambda_{1}>0$. In such a case $\lambda_1$ is the biggest lower frame bound and $\lambda_n$ is the smallest upper frame bound. Moreover $\lambda_{n}=\Vert JJ^*\Vert=\Vert J\Vert^2$. \par
 Now suppose that $[J]_n>0$. By the Cauchy-Binet  theorem $[J]_n=\det(JJ^*)=\prod_{j=1}^n \lambda_j$. Thus
 $$
\lambda_{1}=\frac{\det(JJ^*)}{\prod_{k=2}^n \lambda_k}\ge \frac{[J]_n}{\lambda_{n}^{n-1}}=\frac{[J]_n}{\Vert J\Vert^{2(n-1)}}>0
$$
and $\{v_1,\ldots,v_m\}$ is a frame for $\C^n$ with constants  $\lambda_1\ge {[J]_n}{\Vert J\Vert^{2(1-n)}}$ and \break $\lambda_n=\Vert J\Vert^2$.
 \par
 Conversely, suppose that  $\{v_1,\ldots,v_m\}$ is a frame for $\C^n$ with constants $A$ and $B$. Then $A\le\lambda_1\le\lambda_n\le B$. Hence
 $\Vert J\Vert\le B^{1/2}$ and
 $$
[J]_n=\det(JJ^*)\ge \lambda_1^n\ge A^n>0\hskip.05em  .
$$
This concludes the proof of part (i) of the Lemma.\par
To prove the second part it suffices to observe that  $\{v_1,\ldots,v_n\}$ is a Riesz basis of $\C^n$ if and only if $T$ is an isomorphism and that, in such a case, Riesz constants are also frame bounds. Moreover $[J]_n=\vert\det J\vert^2$ when $m=n$.
\end{proof}
   \begin{theorem}\label{maingen} Suppose that $\frac{\o}{\ell}\le h< \frac{\o}{\ell-\frac{1}{2}}.$ Let $\Phi=\{\varphi_j, 1\le j\le2\ell\} $ be a subset of $\BL$.    Then $\Efi$ is a frame for $\BL$ if and only if  there exist positive constants $\delta,\gamma,\sigma$ and $\eta$ 
such that 
    \begin{equation}\label{maingen:1}
 \delta\le \sum_{j=1}^{2\ell}|\what{\varphi}_j|^2\le\gamma
 \hskip1truecm a.\,e.\  in  \quad (-\o,\o),
\end{equation}
  \begin{equation}\label{maingen:2}
  \Apen{2\ell-1} \ge \sigma  \hskip1truecm a.e.\  in \quad  I_{-}\cup I_{+}\hskip 0.04em, 
\end{equation}
   \begin{equation}\label{maingen:3}
  |\det\,{\Jfi}| \ge \eta  \hskip1truecm a.e.\  in \quad  I.\end{equation}
If   $ h= \frac{\o}{\ell}$ the intervals $I_{-}$ and $I_{+}$ are empty. In this case  $\Efi$ is a Riesz basis  for $\BL$ if and only if    conditions (\ref{maingen:1})   and 
 (\ref{maingen:3}) hold.  
\end{theorem}
\begin{proof}
\par First we shall prove the theorem for   $  \frac{\o}{\ell} <h< \frac{\o}{\ell-1/2}.$ \par Let $\Efi$ be a frame for $\BL$ with frame constants $A$ and $B$. This implies in particular that   $\BL$ coincides with the space $ \Sfi$ generated by $\Efi$. Condition (\ref{maingen:1})  follows from    Proposition \ref{Prop:nec}. Thus we only need to prove  (\ref{maingen:2}) and (\ref{maingen:3}). \par
We recall that the columns of the pre-Gramian $\Jfi$ are the vectors $\cL_h\widehat{\varphi_j}$,\break  $j=1,\ldots,2\ell$ by (\ref{colonnepreg}).
Thus, by Theorem \ref{t: Bownik}(i), the columns of  $\Jfi(x)$ are a frame with constants $A, B$ for the space $\cR_h(B(\omega))(x)$ for a.e. $x$. By Theorem \ref{noteTeoX} we may identify canonically $\cR_h(B(\omega))(x)$ with $\C^{2\ell-1}$ for a.e. $x\in I_-\cup I_+$ and with $\C^{2\ell}$ for a.e.  $x$ in $I$.  Thus, by applying Lemma \ref{simple lemma}(i) with $v_j=\cL_h\widehat{\varphi_j}(x)$ and $J=\Jfi(x)$, we obtain that 
$
[\Jfi(x)]_{2\ell-1}\ge A^{2\ell-1} 
$ for a.e. $x$ in $I_-\cup I_+$ and $\det(\Jfi(x))=[\Jfi(x)]^{1/2}_{2\ell}\ge A^{\ell}$ for a.e. $x$ in $I$. This proves that conditions (\ref{maingen:1})-(\ref{maingen:3}) are necessary.\par
To prove sufficiency  assume that  conditions (\ref{maingen:1})-(\ref{maingen:3}) are  satisfied. 
First we prove that the space $\Sfi$ spanned by $\Efi$ is $\BL$. Since both are $t_o$-shift invariant spaces 
it is enough to show that their range functions coincide almost everywhere.
We recall that the range   $\cR_h(\Sfi) $ of $\Sfi$ is the space spanned by the columns of $\Jfi$. If $x$ is in $I_-$ then by (\ref{colonne:0})  $\cR_h(\Sfi)(x)\subseteq {\rm span }\{\delta_j:\    |j| \le \ell-1  \}$  and the latter space coincides with $\cR_h(\BL)(x)$ by (\ref{RangeterzoX}). On the other hand, $\rango{\Jfi}( x)=2\ell-1$ by (\ref{colonne:0}) and assumption (\ref{maingen:2}). Thus $\cR_h(\Sfi)(x)=\cR_h(\BL)(x)$ because   both have dimension $2\ell-1$. Similar  arguments show  that the range functions coincide almost everywhere also in $I_+$ and in $I$. \par
Next we observe that $\Vert\Jfi(x)\Vert\le \sqrt{2\ell\gamma}$ for a.e. $x$  in $[0,h]$ by (\ref{maingen:1}). Moreover $[\Jfi(x)]_{2\ell-1}>0$ for a.e. $x$ in $I_-\cup I_+$ by (\ref{maingen:2}) and $[\Jfi(x)]_{2\ell}=\vert\det(\Jfi)\vert^2>0$ for a.e. $x$ in $I$ by (\ref{maingen:3}). Thus, by Lemma \ref{simple lemma}(i),  the family $\{\cL_h\widehat{\varphi_1}(x),\ldots, \cL_h\widehat{\varphi_{2\ell}}(x)\}$ is a frame for $\cR_h(B_\omega)(x)$ for a.e. $x$ in $[0,h]$. The upper frame constant  $B=\Vert \Jfi(x)\Vert^2$  is bounded  from  above by $2\ell\gamma$  almost everywhere in $[0,h]$. The lower frame constant $A$ is bounded  from  below by 
$$
[\Jfi(x)]_{2\ell-1}\,\Vert\Jfi(x)\Vert^{2(1-2\ell)}\ge \sigma\,(2\ell\gamma)^{(1-2\ell)}  \quad{\rm a.e.\  in\  }I_-\cup I_+
$$ 
and by 
$$
\vert\det\Jfi(x)\vert^2\,\Vert\Jfi(x)\Vert^{2(1-2\ell)}\ge \eta^2\,(2\ell\gamma)^{(1-2\ell)} \quad{\rm a.e.\  in\ } I.
$$
 Thus $\Efi$ is a frame for $\BL$ by Theorem \ref{t: Bownik}(i).
This concludes the proof of the theorem when $  \frac{\o}{\ell} <h< \frac{\o}{\ell-1/2}.$\par
To prove that if $ h= \frac{\o}{\ell}$ then $\Efi$ is a Riesz basis  for $\BL$ if and only if  
 conditions (\ref{maingen:1})   and 
 (\ref{maingen:3}) hold, one argues in a similar way using Theorem \ref{t: Bownik}(ii) and Lemma \ref{simple lemma}(ii). We omit the details. 
\end{proof}
{\bf Remark} If $\ell=1$ condition (\ref{maingen:2}) is superfluous. Indeed, if $h=\o$  the intervals $I_-$ and $I_+$ are empty. If $\o<h<2\o$ then (\ref{maingen:2}) follows from (\ref{maingen:1}) because
$$[\Jfi]_{2\ell-1}=\sum_{j=1}^2[\what{\varphi_j},\what{\varphi_j}]=
 \sum_{k}\tau_{kh}\sum_{j=1}^2|\what{\varphi_j}|^2
  =\begin{cases} \sum_{j=1}^2|\what{\varphi_j}|^2 &  {\rm  in}\  I_-\\
  {\ }\\
 \tau_{-h}\sum_{j=1}^2 |\what{\varphi_j}|^2  &  {\rm in} \  I_+,\\
 \end{cases}
  $$
and the conclusion follows because   $I_-\subset (-\o,\o)$ and $I_+\subset\tau_{h}(-\o,\o)$.\par
Next we consider the case $  \frac{\o}{\ell-\frac{1}{2}}\le h< \frac{\o}{\ell-1}.$ Then   $\leng{\BL}{t_o}=2\ell-1$ by Corollary \ref{lenBL}. 
Let      $\Phi=$ $\{\varphi_j:1\le j\le 2\ell-1\}$ be a subset of $\BL$  of cardinality $2\ell-1$. By Lemma \ref{pretraslate}  all the rows of the matrix $\Jfi,$ except possibly $\big(\tau_{j h}\what{\varphi}_1,\tau_{j h}\what{\varphi}_2,\dots,\tau_{j h}\what{\varphi}_{2\ell}\big)$, $-\ell\le j\le\ell-1$ vanish. 
 Thus we identify the infinite matrices $\Jfi$,  $\Jfistar$ and $\Gfi$ with their $2\ell-1\times2\ell-1$ submatrices corresponding to these rows. 
The $i$-th column  of $\Jfi,$   $1\le i\le  2 \ell-1$ is 
\begin{equation}\label{colonne:1}
\sqrt{h}\begin{bmatrix} 
\tau_{-(\ell-1)h}\what{\varphi_i}\\
\vdots \\
\what{\varphi_i}\\
\\ \vdots \\
\tau_{(\ell-2)h}\what{\varphi_i}\\
 \tau_{(\ell-1)h}\what{\varphi_i}\\
\end{bmatrix}{\rm{in}}\  K_{-}\hskip0.9truecm
\sqrt{h}\begin{bmatrix} 
\tau_{-(\ell-1)h}\what{\varphi_i}\\
\vdots \\
\what{\varphi_i}\\
\\ \vdots \\
\tau_{(\ell-2)h}\what{\varphi_i}\\
 0\\
\end{bmatrix}\ {\rm{in}}\  K\hskip0.9truecm
\sqrt{h}\begin{bmatrix} \tau_{-\ell h}\what{\varphi_i}\\
\tau_{-(\ell-1)h}\what{\varphi_i}\\
\vdots \\
\what{\varphi_i}\\
\\ \vdots \\
\tau_{(\ell-2)h}\what{\varphi_i} \end{bmatrix}\  {\rm{in}}\  K_{+}.
\end{equation}
\par
\begin{theorem}\label{maingendue} Let $\Phi \subset \BL$ such that $\Phi=\{\varphi_j, 1\le j<2\ell-1\},$ and $\ell\not=1$ such that 
$\frac{\o}{\ell-\frac{1}{2}}\le h< \frac{\o}{\ell-1}. $
Then $\Efi$ is a frame for $\BL$ if and only if  there exist positive constants $\delta, \gamma, \sigma$ and $\eta$ such that
    \begin{equation}\label{maingendue:1}
\delta\le \sum_{j=1}^{2\ell-1}|\what{\varphi}_j|^2\le\gamma
 \hskip1truecm a.e. \   in\  (-\o,\o),
\end{equation}
  \begin{equation}\label{maingendue:2}
\Apen{2\ell-2}\ge \sigma  \hskip1truecm a.e.\    in \ K, \end{equation}
   \begin{equation}\label{maingendue:3}
|\det\, \Jfi|\ge \eta  \hskip1truecm a.e.\   in \   K_{-}\cup K_{+}.
\end{equation}
If $h=\frac{\o}{\ell-\frac{1}{2}}$ then $\Efi$ is a Riesz basis if and only if conditions  (\ref{maingendue:1}) and (\ref{maingendue:3}) hold. 
\end{theorem}
The proof is similar to that of Theorem \ref{maingen}. We omit the details.
  \section{The dual generators}
  \label{s: Dual generators} 
 Let $N=\leng{\BL}{t_o}$ be the length of $\BL$ as $t_o$-shift-invariant space  and let $\Phi=\{\varphi_1,\ldots,\varphi_N\}$ be a subset of $\BL$.   In this section we shall  find   the dual generators  $\Phi^*$ when $\Efi$ is a Riesz basis or a frame of $\BL$.
With a slight abuse of notation in this section we shall denote by ${\Phi}$  the vector $({\varphi_1},\ldots,{\varphi_N})$ and by ${\Phi^*}$ the vector $({\varphi^*_1},\ldots{\varphi^*_N})$.
\par
It is well known that if  $\Efi$ is a Riesz basis for $\Sfi$  then the  Gramian matrix is invertible and  the Fourier transform of the dual generators are given by  
\begin{equation}\label{dualRiesz}
\what{\Phi^*}^\top= \overline{\Gfimu}{\what{\Phi}}^\top,
\end{equation} 
where $ v^{\top}$ denotes  the transpose of the vector $v$. In  Theorems \ref{maindualnew} - \ref{maindual4gen} 
 we give explicit formulas for  the Fourier transforms of  the dual generators  when
   $\Efi $ is a frame satisfying the hypothesis of Theorems \ref{maingen} or  \ref{maingendue}.   The proof is based on   the dual Gramian matrix $\Gfitilde$ representation  of the operator  $\Tfi \Tfia. $  From (\ref{LFTTsF})   we obtain 
\begin{equation}\nonumber
\cL_h  { \mathcal F} \Tfi \Tfia { \mathcal F}^{-1} \what{\varphi_k}^*= \Jfi \Jfistar  \  \cL_h \what{\varphi_k}^* \hskip 1truecm k=1,\ldots,N.
\end{equation}
By (\ref{fiduale})   the left hand side  is equal to 
$\cL_h   \what{\varphi_k}.$ Hence   $$ \cL_h   \what{\varphi_k} = \Jfi \Jfistar  \  \cL_h   \what{\varphi_k^*},\hskip 1truecm k=1,\ldots.N
 $$
which, by (\ref{colonnepreg}),
can be written
 \begin{equation}\label{generators:1}
 \Jfi=\Jfi\Jfistar\, \Jfidual.
 \end{equation}

As in Section 2 we identify the infinite matrices $\Jfi$, $\Jfistar$ and $\Jfidual$ with $N\times N$ matrices by neglecting their vanishing  rows and columns (see the discussion after Lemma \ref{pretraslate}). Thus we shall interpret (\ref{generators:1})  as an identity between $N\times N$ matrices.
\par  Under the assumptions of Theorems \ref{maingen}  and \ref{maingendue} the interval $[0,h]$ is the disjoint union of three intervals where the pre-Gramian is either invertible or has rank $N-1$. In the latter case either the first or the last row of the pre-Gramian vanishes. In this case we shall denote by $\Jbfi$  the $(N-1)\times N$ submatrix of $\Jfi$ obtained by deleting the vanishing row from $\Jfi$. It is straightforward to see that in this case equation (\ref{generators:1}) reduces to
 \begin{equation}\label{generators:2}
{\mathbb J}_{\Phi,t_o}=\Jbfi \Jbfistar\, \Jbfidual . \end{equation}
  We  regard (\ref{generators:1}) and (\ref{generators:2}) as  equations for the unknowns   $\Jfidual$ and ${\mathbb J}_{\Phi^*,t_o}$, respectively.  
In the intervals  where the matrix $ \Jfi$ is invertible we can solve for $\Jfidual$ in (\ref{generators:1}), obtaining  that 
\begin{equation}\label{Gei:1}
\Jfidual=(\Jfistar)^{-1}.\end{equation}  In the intervals where the rank of $\Jfi$ is $N-1$  we can solve for the submatrix
$\Jbfidual$   obtaining that
\begin{equation}\label{Gei:2}
\Jbfidual=(\Jbfi \Jbfistar )^{-1}\Jbfi.
\end{equation} 
We  recall that if $A$ is a $N\times(N-1)$ matrix of rank $N-1$ then  its  Moore-Penrose inverse $A^\dagger$    is 
 \begin{equation}\label{defpenrose}
 A^\dagger=(A^*A)^{-1} A^*
 \end{equation} 
(see \cite{BIG}).   Therefore by (\ref{Gei:2})
\begin{equation}\label{Gei:3}
 \Jbfidual= (\Jbfistar)^\dagger . 
\end{equation}
 We refer the reader to \cite{BIG} for the definition and the properties of  the Moore-Penrose  inverse of a  matrix.
\par
By using (\ref{Gei:1}) and (\ref{Gei:2}) we shall obtain explicit formulas for the Fourier transforms of the dual generators. For the sake of clarity first we state and prove the result for $N=2,3,4$.  By Corollary \ref{lenBL} these cases correspond to $\o\le h<2\o$, $\frac{\o}{2}\le h<\frac{2}{3}\o$ and $\frac{2}{3}\o \le h<  \o$ respectively.   
   \begin{theorem}\label{maindualnew} Assume that $\o\le h<2\o$ and 
let $\Phi$ denote the vector $(\varphi_1, \varphi_2 )$ where  $\varphi_1, \varphi_2\in\BL$.  If       (\ref{maingen:1}) and   (\ref{maingen:3}) hold with $\ell=1$, i.e. if  $\Efi$ is a frame for $\BL$, then
\begin{equation}\label{sol:2fin}
\what{\Phi^*}=\begin{cases}
\quad D\,\tau_h\ove{\what{\Phi}}^\perp &  { in}\   [-\o,\o-h]\\
\quad{h^{-1}\Vert{\what{\Phi}}\Vert^{-2}}\  {\what{\Phi}}
&  { in}\  (\o-h,h-\o) \\
 -D\,\tau_{-h}\ove{\what{\Phi}}^\perp &  { in}\ [h-\o,\o]
\end{cases}
\end{equation}
where $D=(\det\ {\Jfistar})^{-1}$ and $\what{\Phi}^\perp=(\what{\varphi_2}, -\what{\varphi_1})$.
 Note that  if $h=\o$ the central interval is empty.
\end{theorem}
 \begin{proof} 
 Assume first that $\o< h<2\o.$   We recall that $I_-=(0,h-\o)$, $I=(h-\o,\o)$ and $I_+=(\o,h)$. The pre-Gramian $\Jfi$ is 
 \begin{equation}\nonumber
\sqrt{h}\begin{bmatrix}
0 & 0\\
 \what{\varphi_1} & \what{\varphi_2}   \\
\end{bmatrix}{\rm in}\hskip.3truecm  I_-
\hskip1.4truecm    
 \sqrt{h} \begin{bmatrix}
 \tau_{-h}\what{\varphi_1} &  \tau_{-h}\what{\varphi_2}\\
0 & 0  \\
\end{bmatrix}\ {\rm in} \hskip.3truecm I_+\end{equation}
 \begin{equation}\nonumber
 \sqrt{h}  \begin{bmatrix}
 \tau_{-h} \what{\varphi_1} & \tau_{-h} \what{\varphi_2}\\
 \what{\varphi_1} & \what{\varphi_2}  \\ \end{bmatrix}{\rm in}\hskip.2truecm  I.  
 \end{equation}
 The same formulas hold for $\Jfidual$ with $\what{\varphi_i}$ replaced by $\what{\varphi^*_i}, j=1,2$.
Therefore  
\begin{equation}\label{Gei:4}
\Jbfi=\begin{cases}
\sqrt{h}\, \what{\Phi} & {\rm in\ } I_-  \\
\sqrt{h}\, \tau_{-h}\what{\Phi} & {\rm in\ } I_+ \\
\end{cases}
\hskip1truecm
\Jbfidual=\begin{cases}
\sqrt{h}\, \what{\Phi^*} & {\rm in\ } I_-  \\
\sqrt{h}\, \tau_{-h}\what{\Phi^*} & {\rm in\ } I_+\,. \\
\end{cases}
\end{equation} 
  By  assumptions   (\ref{maingen:1}) and  (\ref{maingen:3})
   ${\Jfi}$ has rank $1$   in $I_{-}\cup I_{+}$ and  rank $2$ in  $I.$   Hence $\Jbfidual=(\Jbfistar)^\dagger$ in $I_-\cup I_{+}$ and  
    \begin{equation}\label{Gei:6}
 \Jfidual={\Jfistar}^{-1} \qquad  {\rm \ in} \quad I.  \end{equation}
  First we find $\Jbfidual$ in $I_-\cup I_+$. By using (\ref{Gei:2}) we get
    \begin{equation}\nonumber
\Jbfidual=(\Jbfistar)^\dagger   =\begin{cases}
\frac{\what{\Phi}}{\sqrt{h}\Vert{\what{\Phi}}\Vert^2}
 & {\rm in\ } I_-  \\
 \\
\frac{\tau_{-h}\what{\Phi}}{\sqrt{h}\Vert{\tau_{-h}\what{\Phi}}\Vert^2}
 & {\rm in\ } I_+ \hskip0.05em .\\
\end{cases}
\end{equation}
By (\ref{Gei:4}) we obtain that
$\what{\Phi^*}= 
\frac{\what{\Phi}}{h\Vert{\what{\Phi}}\Vert^2}$ in $I_- $ and 
$\tau_{-h}\what{\Phi^*}= \frac{\tau_{-h}\what{\Phi}}{h\Vert{\tau_{-h}\what{\Phi}}\Vert^2} $ in $   I_+$.  Since $\tau_{-h}I_+=(\o-h,0)$
we finally get 
  \begin{equation}\nonumber
\what{\Phi^*}=    
\frac{\what{\Phi}}{h\Vert{\what{\Phi}}\Vert^2}
 \hskip1.5truecm  {\rm in\ } (\o-h,h-\o).   \end{equation}
Next we find the dual generators in the remaining intervals. By     (\ref{Gei:6})  
\begin{equation}\nonumber
\Jfidual=\sqrt{h}\ ({\det} {\Jfistar})^{-1} \begin{bmatrix}
\ove{\what{\varphi_2}} & -\ove{\what{\varphi_1}}    \\
-\tau_{-h}\ove{\what{\varphi_2}}   & \tau_{-h}\ove{\what{\varphi_1}}   \\
\end{bmatrix} \hskip1truecm {\rm a. e. \ in}\hskip.3truecm (h-\o,\o).
\end{equation} 
By translating and reminding that the  pre-Gramian matrix is $h$-periodic
\begin{equation}\nonumber
\what{\Phi}^*=\begin{cases}
{D}\,\tau_{h}\big(\ove{\what{\varphi_2}},
-\ove{\what{\varphi_1}}\big)={ D}\,\tau_{h}\ove{\what{\Phi}}^{\,\perp}&\quad {\rm in}\quad (-\o,\o-h)\\
 D\,\tau_{-h}\big(-\ove{\what{\varphi_2}},
\ove{\what{\varphi_1}} \big)= -{ D}\,\tau_{-h}\ove{\what{\Phi}}^{\,\perp}&\quad {\rm in}\quad (h-\o,\o).
\end{cases}\end{equation}
 This completes the proof of the theorem when $\o< h< 2\o$.
\par If $h=\o$  one argues as before; the only difference is that now the interval $(h-\o,\o-h)$ is empty. \end{proof}

To find an explicit expression of the dual generators when $N>2$ we need  formulas for the rows of the Moore-Penrose  inverse of a $N\times (N-1)$ matrix of full rank.
 Given  a $(N-1)$-ple of  vectors $ (U_1,\ldots,U_{N-1})$ in  $\C^N$     their cross  product is 
 $$\crossp{1}{N-1} U_j=U_1 \cp U_{2}\dots\cp U_{N-1} =
 \det \begin{bmatrix}
  \mathbf{e}_1\quad &   \mathbf{e}_2 \quad&{ \ }  &{\dots}  &  \mathbf{e}_N\quad \\
U_1^1\quad& U_1^2 \quad &{ \ }      &{\dots}   & U_1^N\quad \\
\vdots &  \vdots &{ \ }    &  &  \vdots \\
U_{N-1}^1 & U_{N-1}^2 &{ \ }       &{\dots}   & U_{N-1}^N  \\ \end{bmatrix} 
$$
 where the $\{\mathbf{e_j}:j=1,\ldots,N\}$ is the canonical basis of $\C^N$.
 Notice that if $N=2$ then $\crossp{}{}U=U^\perp$.   Given   a vector  $W $ in $\C^N $ and an integer $k\in\{1,2,\ldots, N-1\}$,\break   we shall denote by $ \displaystyle{\crosspw{j=1}{N-1}{U_j}{U_k}{W}} $ the cross product of the $(N-1)$-ple \break $(U_1,\ldots, U_{k-1},W,U_{k+1},\ldots,U_{N-1})$, i.e. 
 $$ \crosspw{j=1}{N-1}{U_j}{U_k}{W}  = U_1 \cp U_{2}\dots\cp U_{k-1}\cp W\cp U_{k+1}\dots \cp U_{N-1}.$$
  \begin{lemma} \label{inversa}
   Le $M$ be an $n\times n $ invertible matrix. Denote by $R_j, \, j=1,\dots,n$, its rows and by  $ C_j $ its   columns. Then  the columns of $M^{-1}$ are  $$(-1)^{k+1}({\det} M)^{-1}\crosspd{j=1}{n} {j\not = k}R_j  \hskip1truecm 1\le k\le n $$
 and the rows are 
 $$(-1)^{k+1}( {\det} M)^{-1} \crosspd{j=1}{n} {j\not = k}C_j  \hskip1truecm 1\le k\le n. $$
  \end{lemma}
  \begin{proof}  Let $M_{kj}$ be the $kj$-cofactor of $M$. Then the $k$-th row of the matrix ${\rm cof}(M)$  of cofactors of $M$ is
\begin{align*}
 L_k&=(-1)^{k+1}(M_{k1}{\mathbf e}_1-M_{k2}{\mathbf e}_2+\cdots+(-1)^{n-1}M_{kn}{\mathbf e}_n) \\ 
&= (-1)^{k+1}\crosspd{j=1}{n} {j\not = k}R_j.
\end{align*}
The first identity follows from the fact that $M^{-1}=(\det M)^{-1}\,{\rm cof} (M^{\top})$. The proof of the second identity is similar.
\end{proof}
  \begin{lemma}\label{penrose}
   Le $A$ be a  $n\times (n-1)$  complex  matrix of maximum rank. Denote by $A_j,$ $ j=1,\dots,n-1,$ the columns of $A$  and by  $ P_k, \, k=1,\dots,n-1$   the  rows of $A^\dagger.$   Then  \begin{equation}\label{penrose:1} P_k= (-1)^{n}[{\det} (A^*A)]^{-1}\crosspw{j=1}{n-1}  {A_j}{A_k}{W}  \hskip1truecm 1\le k\le n-1 
   \end{equation}
   where $\displaystyle{W=\crossp{j=1}{n-1}\overline{A}_j}$.
    \end{lemma}
   \begin{proof}  Since $\rango{A}=n-1$ the null space of $A^*$ is the space spanned by $W$. Let $A_b$ the matrix obtained by bordering $A$ with the column $W$, i.e.
   $$
A_b=\left[\begin{array}{cc}A  &W \end{array}\right].
$$ 
Then $A_b$ is invertible because ${\det} A_b=\displaystyle{\crossp{j=1}{n-1}} A_j \cdot W =|W|^2$. By  \cite[Thm. 8]{BIG}
$$
A_b^{-1}=\left[\begin{array}{ll}A^\dagger  \\
W^\dagger \end{array}\right].
$$ 
Thus, for every $k=1,\ldots,n-1$, the $k$-th row of $A^\dagger$ is the $k$-th row of $A_b^{-1}$. 
Hence, by Lemma \ref{inversa} and 
the anticommutativity   of the cross product, we obtain that   
\begin{align*}\label{penrose:2}
 P_k=& (-1)^{k+1} \big({\det}A_b\big)^{-1}  A_1\cp A_2 \cp\dots \cp A_{k-1}\cp A_{k+1}\dots \cp A_{n-1} \cp W  \\ 
=& (-1)^n |W|^{-2}\crosspw {j=1}{n-1}{A_j}{A_k}{W}.  
\end{align*} 
To conclude the proof we observe that $|W|^2$ is the sum of the squares of the absolute values of the minors of   order  $n-1$ of $A$, i.e. the determinant of $ A^*A$, by the Cauchy-Binet formula.
 \end{proof}
 \begin{theorem}\label{maindual3}
Assume that $\frac{2}{3}\o \le h <\o $ and let $\Phi$ denote the vector $(\varphi_1,\varphi_2,\varphi_3)$ where $\varphi_j\in \BL$, $j=1,2,3.$ If  assumptions (\ref{maingendue:1})-(\ref{maingendue:3}) hold with $\ell=2$, i.e. if $\Efi$is a frame for $\BL$,  then
 the Fourier transform of the   dual generators \break $\Phi^*=(\varphi_1^*,\varphi_2^*,\varphi_3^*)$ is 
 \begin{equation}\label{maindual3t}
\widehat{\Phi}^*= \begin{cases} 
  D\quad   \tau_{h} \overline{\what\Phi}\cp \tau_{2h} \overline{\what\Phi} &{\rm{in}}\ (-\o,\o-2h)\\
   \Du  \  \big( \tau_{h} {\what\Phi}\cp 
  {\what\Phi}\big)\cp \tau_{h}\overline{\what\Phi} &{\rm{in}} \ (\o-2h,h-\o)\\
D  \quad \tau_h \overline{\what\Phi}\cp 
  \tau_{-h} \overline{\what\Phi} &{ \rm{in}}\  (h-\o,\o-h)\\
-  \Du  \    \tau_{-h} \overline{\what\Phi} \cp  \big( \tau_{-h} {\what\Phi}\cp 
 {\what\Phi}\big) &{\rm{in}} \ (\o-h,2h-\o)\\
 D\quad   \tau_{-2h} \overline{\what\Phi}\cp   \tau_{-h} \overline{\what\Phi}&{\rm{in}}\  (2h-\o,\o)\\
  \end{cases}
  \end{equation}
   where $\Du =h\,\JJbstinv\ $ and   $D=\sqrt{h}\,({\det}\ {\Jfistar})^{-1}$.  Note that  if $h=\frac{2\o}{3}$ the    intervals  $(\o-h,2h-\o)$  $(\o-2h,h-\o)$ are empty. \end{theorem}
 \begin{proof}
Assume first that  $\frac{2\o}{3}<h<  {\o}. $ We recall  that    $K_{-}, K$ and $ K_{+}$   denote  the intervals   $(0,\o-h)$, $(\o-h,2h-\o)$ and $(2h-\o,h)$ respectively, defined in (\ref{K}).  By (\ref{colonne:1})  the  matrix    $\Jfi$ is 
\begin{equation}\label{maindual3:0}
\sqrt{h}\begin{bmatrix}\tau_{-h}\what{\Phi}\\
\what{\Phi}\\
\tau_{h}\what{\Phi}
\end{bmatrix}{\rm{in}}\  K_{-}\hskip0.9truecm
\sqrt{h}\begin{bmatrix}\tau_{-h}\what{\Phi}\\
 \what{\Phi}\\ 0\\
\end{bmatrix}\ {\rm{in}}\  K\hskip0.9truecm
\sqrt{h}\begin{bmatrix}\tau_{-2h}\what{\Phi}\\
\tau_{-h}\what{\Phi}\\
\what{\Phi}
\end{bmatrix}\  {\rm{in}}\  K_{+}.
\end{equation}
The same formulas hold for $\Jfidual$ with $\what{\Phi}$ replaced by $\what{\Phi^*} $.
Therefore  
\begin{equation}\label{maindual3:1}
\Jbfi=\sqrt{h}\begin{bmatrix}\tau_{-h}\what{\Phi}\\
 \what{\Phi}\\  
\end{bmatrix} \hskip0.3truecm {\rm and }\hskip0.3truecm
\Jbfidual =\sqrt{h}\begin{bmatrix} 
\tau_{-h}\what{\Phi^*}\\
\what{\Phi^*}
\end{bmatrix}\  {\rm{in}}\  K.
\end{equation}
 By assumptions 
(\ref{maingendue:2}) and (\ref{maingendue:3})  the matrix $\Jfi$ has   rank   $3$  in $K_-\cup K_+ $ and rank 2 in $K$. Hence 
$ \Jfidual={\Jfistar}^{-1}$ in $ K_-\cup K_+ $
 and $\Jbfidual=(\Jbfistar)^\dagger$ in $K$. 
  \par First we find $\Jfidual$ in $K_-$. By  Lemma \ref{inversa}   $$ \tau_{-h}\what{\Phi^*}=
  D\     \overline{\what\Phi} \cp   \tau_{h}\overline{\what\Phi}, \hskip.45 truecm
 \what{\Phi^*}=
  D\   \tau_h \overline{\what\Phi}\cp  \tau_{-h} \overline{\what\Phi},  \hskip.45 truecm
  \tau_{h}\what{\Phi^*}=
  D\   \tau_{-h}\overline{\what\Phi}\cp  
 \overline{\what\Phi} \hskip.5truecm  {\rm a.e. \ in}\   K_- 
$$
where 
$D=\sqrt{h} \,({\det}\ {\Jfistar})^{-1}$.
By translating  the first and the last  identities and reminding that the pre-Gramian is  $h$-periodic, we obtain 
\begin{equation}\label{maindual3:3}  \what{\Phi^*}=\begin{cases}
 D \hskip.4truecm    \tau_{h}
 \overline{\what\Phi} \cp    \tau_{2h}\overline{\what\Phi}\hskip .2truecm 
&{\rm {in}}\  (-h,\o-2h)\\
  D\quad \,  \tau_h \overline{\what\Phi}\cp 
  \tau_{-h} \overline{\what\Phi}
  &{\rm{in}} \ (0,\o-h)\\
 D\    \tau_{-2h}
 \overline{\what\Phi} \cp     \tau_{-h}\overline{\what\Phi}\hskip 0.2truecm 
&{\rm {in}}\ (h,\o).\\
\end{cases}\end{equation}
The same calculation  in $K_{+}$ gives
$$ \tau_{-2h}\what{\Phi^*}=
    D\     \tau_{-h} \overline{\what\Phi} \cp  
   \overline{\what\Phi} \hskip1truecm
 \tau_{-h}\what{\Phi^*}=
   D\    \overline{\what\Phi}\cp 
  \tau_{-2h} \overline{\what\Phi} \hskip1truecm
 \what{\Phi^*}=
    D\   \tau_{-2h}\overline{\what\Phi}\cp  
 \tau_{-h} \overline{\what\Phi}.
$$ By translating the first two identities we obtain
\begin{equation}\label{maindual3:4}  \what{\Phi^*}=\begin{cases}
    D\hskip.4truecm    \tau_{h}
 \overline{\what\Phi} \cp    \tau_{2h}\overline{\what\Phi}\hskip .2truecm 
&{\rm {in}}\ (-\o,-h)\\
   D\quad \,  \tau_h \overline{\what\Phi}\cp  
  \tau_{-h} \overline{\what\Phi}
  &{\rm{in}} \ (h-\o,0)\\
    D\    \tau_{-2h}
 \overline{\what\Phi} \cp      \tau_{-h}\overline{\what\Phi}\hskip 0.2truecm 
&{\rm {in}}\ (2h-\o,h).\\
\end{cases}\end{equation}
Next we  find $\Jbfidual$ in $K$. By   (\ref{maindual3:1}) and   Lemma \ref{penrose} the  rows of Moore-Penrose inverse of   $\Jbfistar$   are  
$$   
- \sqrt{h} \Du \   \big(\tau_{-h}  {\what\Phi} \cp 
  {\what\Phi}\big)\cp    \overline{\what{\Phi}} \hskip2truecm
 -\sqrt{h} \Du \   \tau_{-h} \overline{\what\Phi} \cp  \big(\tau_{-h}{\what\Phi}\cp  {\what\Phi} \big) 
$$
where $\Du =h\, \JJbstinv $. Hence, by  (\ref{maindual3:1}),  $$   
\tau_{-h}\what{\Phi^*}= - \Du\    \big(\tau_{-h}  {\what\Phi} \cp 
  {\what\Phi}\big)\cp  \overline{\what{\Phi}} \hskip1truecm
 \what{\Phi^*}=- \Du\    \tau_{-h} \overline{\what\Phi} \cp  \big(\tau_{-h}{\what\Phi}   \cp   
  {\what\Phi} \big).
$$
By translating the first identity and using the anticommutativity   of the cross pro\-duct,  we obtain \begin{equation}\label{maindual3:5}  \what{\Phi^*}=\begin{cases}
  \, \Du\  \tau_{h}  \overline{\what{\Phi}} \cp \big(   {\what\Phi} \cp  \tau_{h}  {\what\Phi}\big)   \hskip .2truecm 
&{\rm {in}} \ ( \o-2h,h-\o)\\
 \, \Du\ \big(\tau_{-h}{\what\Phi}\cp  
  {\what\Phi} \big)\cp \tau_{-h} \overline{\what\Phi} 
  &{\rm{in}}\  (\o-h,2h-\o).\\
\end{cases}\end{equation}
The conclusion follows from formulas (\ref{maindual3:3}), (\ref{maindual3:4}) and (\ref{maindual3:5}).
 This completes the proof of the theorem when $\frac{2\o}{3}< h< \o$.
\par If $h=\frac{2\o}{3}$  one argues as before; the only difference is that now the interval $(\o-h,2h-\o)$ is empty.   \end{proof}
\begin{theorem}\label{maindual4}
Assume that $\frac{\o}{2}\le h< \frac{2}{3}\o $ and let $\Phi$ denote the vector $(\varphi_1,\varphi_2,\varphi_3,\varphi_4)$ where $\varphi_j\in \BL$, $j=1,\dots,4.$ If  (\ref{maingen:1})-(\ref{maingen:3}) hold with $\ell=2$, i.e. if $\Efi$  is a frame for $\BL$,
then the Fourier transform of the   dual generators $\Phi^*=(\varphi_1^*, \varphi_2^*, \varphi_3^*,\varphi_4^*)$  is   \begin{equation}\label{maindual4:00}
\widehat{\Phi}^*= \begin{cases} 
  \Du   \quad \tau_{-h} \overline{\what\Phi} \cp  
\big( \tau_{-h} {\what\Phi}   \cp \what{\Phi} \cp  \tau_{h} {\what\Phi} \big)\cp   \tau_h \overline{\what{\Phi} } &{  {in}}\  (0,2h-\o)\\
  D\quad   \tau_{-2h} \overline{\what\Phi} \cp 
  \tau_{-h} \overline{\what\Phi}  \cp \tau_{h}\overline{\what\Phi} &{ {in}} \ (2h-\o,\o-h)\\
 \Du   \quad     \tau_{-2h} \overline{\what\Phi} \cp 
  \tau_{-h} \overline{\what\Phi}  \cp \big(\tau_{-2h} {\what\Phi} \cp   \tau_{-h }\what{\Phi} \cp   {\what\Phi} \big)&{ {in}} \ (\o-h,3h-\o)\\
  D\quad      \tau_{-3h} \overline{\what\Phi}   \cp  \tau_{-2h} \overline{\what\Phi} \cp \tau_{-h}  \overline{\what\Phi} &{{in}}\  (3h-\o,\o)\\
  \end{cases}
  \end{equation}
  where $\Du=h^2\,\JJbstinv$ and  $D= {h}\,({\det}\overline{\Jfi})^{-1}.$ 
 The expression of $\Phi^*$ in    $(-\omega,0)$ is obtained by reflecting each interval in (\ref{maindual4:00}) around zero and   replacing  $\tau_{jh}$ with $\tau_{-jh}$ in the expression of  $\what\Phi^*$ in the reflected interval. Note that if $h=\frac{\o}{2}$ then the first and third  intervals are empty.
  \end{theorem}
 \begin{proof} Assume first that  $\frac{\o}{2}< h< \frac{2\o}{3}.$ We recall that by $I_-,I $ and $I_+$ we denote the intervals  $  (0,2h-\o)$, $(2h-\o,\o-h)$ and   $(\o-h,h) $  defined in (\ref{I}). 
By (\ref{colonne:0}) the  columns of $\Jfi$  are 
\begin{equation}\nonumber
\sqrt{h}\begin{bmatrix} 0\\ \tau_{-h}\what{ \Phi}\\
\what{ \Phi}\\
\tau_{h}\what{ \Phi}
\end{bmatrix}{\rm{in}}\  I_{-}\hskip0.9truecm
\sqrt{h}\begin{bmatrix}\tau_{-2h}\what{ \Phi}\\
\tau_{-h}\what{ \Phi}\\
  \what{ \Phi}\\
  \tau_{h}\what{ \Phi}\\
\end{bmatrix}\ {\rm{in}}\  I\hskip0.9truecm
\sqrt{h}\begin{bmatrix}\tau_{-2h}\what{ \Phi}\\
\tau_{-h}\what{ \Phi}\\ 
 \what{ \Phi}\\ 0\\
\end{bmatrix}\  {\rm{in}}\  I_{+}.
\end{equation}
The same formulas hold for $\Jfidual$ with $\what{\Phi}$ replaced by $\what{\Phi^*} $.
Therefore \begin{equation}\label{maindual4:1b}
\Jbfi=\sqrt{h}\begin{bmatrix}\tau_{-h}\what{\Phi}\\
 \what{\Phi}\\ \tau_{h}\what{\Phi} 
\end{bmatrix} \hskip0.6truecm {\rm and } \hskip0.6truecm
\Jbfidual =\sqrt{h}\begin{bmatrix} 
\tau_{-h}\what{\Phi^*}\\
\what{\Phi^*}\\ \tau_{h}\what{\Phi^*}
\end{bmatrix}\quad  {\rm{in}}\ I_
-\end{equation}
 \begin{equation}\label{maindual4:2b} \Jbfi=\sqrt{h}\begin{bmatrix}\tau_{-2h}\what{\Phi}\\
\tau_{-h} \what{\Phi}\\ \what{\Phi} 
\end{bmatrix}\hskip0.6truecm {\rm and } \hskip0.6truecm
\Jbfidual =\sqrt{h}\begin{bmatrix} 
\tau_{-2h}\what{\Phi^*}\\
\tau_{-h}\what{\Phi^*}\\ \what{\Phi^*}
\end{bmatrix}\quad {\rm{in}}\  I_+.
\end{equation}
  By assumptions 
(\ref{maingen:2}) and (\ref{maingen:3})  the matrix $\Jfi$ has   rank   $3$  in  $\, I_-\cup I_+ \,$ and rank 4 in $I$. Hence  
 $\Jbfidual=(\Jbfistar)^\dagger$ in $ I_-\cup I_+  $ and $ \Jfidual={\Jfistar}^{-1}$ in $I$.
\par First we find $\Jfidual$ in $I$;  by Lemma \ref{inversa}
    \begin{align*}
 \tau_{-2h}\what{\Phi^*}=&\quad\, D \quad\,  \tau_{-h}\ove 
 {\what{\Phi}}\cp \ove {\what{\Phi}}\cp\tau_{h}\ove {\what{\Phi}}\\
  \tau_{-h}\what{\Phi^*}=&-\,  D\quad \tau_{-2h}\ove{\what{\Phi}}\cp \ove{\what{\Phi}}\cp\tau_{h}\ove{\what{\Phi}} \\
    \what{\Phi^*}=&\quad\, D\quad\, \tau_{-2h}\ove{\what{\Phi}}\cp  
  \tau_{-h}\ove{\what{\Phi}}\cp\tau_{h}\ove{\what{\Phi}}\\ \tau_{h}\what{\Phi^*}=&-\, D \quad \tau_{-2h}\ove{\what{\Phi}}\cp \tau_{-h}\ove{\what{\Phi}}\cp\ove{\what{\Phi}}\end{align*}
in $(2h-\o,\o-h);$  here   $D= {h}\,({\det}\overline{\Jfi})^{-1}.$ By translating   and reminding that the  matrix $\Jfi$ is $h$-periodic we obtain
    \begin{equation}\label{maindual4:1}
 \what{\Phi^*}=\begin{cases} \quad \, D \quad \tau_{h}\ove 
 {\what{\Phi}}\cp \tau_{2h}\ove {\what{\Phi}}\cp\tau_{3h}\ove {\what{\Phi}}& {\rm{in}}\quad(-\o,\o-3h) \\
 -\, D \quad \tau_{-h}\ove{\what{\Phi}}\cp   \tau_{h}\ove{\what{\Phi}}\cp\tau_{2h}\ove{\what{\Phi}}& {\rm{in}}\quad (h-\o,\o-2h)\\
 \quad\, D \quad \tau_{-2h}\ove{\what{\Phi}}\cp  
  \tau_{-h}\ove{\what{\Phi}}\cp\tau_{h}\ove{\what{\Phi}}& {\rm{in}}\quad(2h-\o,\o-h)\\
 -\, D \quad \tau_{-3h}\ove{\what{\Phi}}\cp \tau_{-2h}\ove{\what{\Phi}}\cp\tau_{-h}\ove{\what{\Phi}}
& {\rm{in}}\quad(3h-\o,\o). 
 \end{cases}
\end{equation} 
Notice that if $[a,b]$ is any of the intervals in the r.h.s. of (\ref{maindual4:1}) the expression of $ \what{\Phi^*}$ in $[a,b]$  can be obtained  from that in $[-b,-a]$, by replacing $h$ by $-h$  in the translations $\tau_{jh}, |j| \le3.$
 \par
Next we find the dual generators in the remaining intervals. First let us consider the interval $I_{-}.$ Here by  Lemma \ref{penrose}    the rows of  
the Moore-Penrose inverse  of  $\Jbfistar$ are   
$$ \sqrt{h} \,\Du\  W\cp\ove{\what\Phi}\cp\tau_{h}\ove{\what\Phi},
\hskip0.7truecm\sqrt{h}\, \Du\   \tau_{-h}\ove{\what\Phi}\cp W \cp\tau_{h}\ove{\what\Phi},\hskip0.7truecm \sqrt{h}\, \Du \  \tau_{-h}\ove{\what\Phi}\cp\ove{\what\Phi}\cp W
$$  where $W=\tau_{-h}{\what\Phi}\cp{\what\Phi} \cp\tau_{h}{\what\Phi}
 $   and $\, \Du=h^2\,\JJbstinv.$ By using (\ref{maindual4:1b})  and  translating we obtain  
 \begin{equation}\label{maindual4:2}
 \what{\Phi^*}=\begin{cases}   
  \Du   \quad \tau_h W\cp\tau_h {\what{\Phi}}\cp  \tau_{2h}\ove{\what{\Phi}} \  \quad  \qquad{\rm{in}}& (-h,h-\o)\\
 \Du \quad\tau_{-h}\ove{\what{\Phi}}\cp  
 W\cp\tau_{h}\ove{\what{\Phi}}\ \  \quad \qquad{\rm{in}}& (0,2h-\o)\\ 
  \Du \  \quad \tau_{-2h}\ove{\what{\Phi}}\cp  \tau_{-h}\ove{\what{\Phi}}\cp \tau_{-h}W \quad{\rm{in}}& (h,3h-\o).\\
   \end{cases}
\end{equation}  
Finally we  consider the interval $I_{+}$. Here    the rows of the Moore-Penrose inverse of 
  $\Jbfistar$ are 
 $$ \sqrt{h} \Du\ W_o\cp\tau_{-h}\ove{\what\Phi}\cp\ove{\what\Phi},
\hskip0.55truecm\sqrt{h}  \Du\  \tau_{-2h}\ove{\what\Phi}\cp W_o\cp\ove{\what\Phi}, \hskip0.50truecm \sqrt{h}  \Du\ \tau_{-2h}\ove{\what\Phi}\cp\tau_{-h}\ove{\what\Phi}\cp W_o
$$ 
where $W_o=\tau_{-2h}{\what\Phi}\cp\tau_{-h}{\what\Phi} \cp{\what\Phi}
.$ By using (\ref{maindual4:2b}) and translating we obtain 
   \begin{equation}\nonumber
 \what{\Phi^*}=\begin{cases} \Du\quad    \tau_{2h}W_o\cp\tau_{h}\ove{{\what{\Phi}}}\cp\tau_{2h} \ove{\what{\Phi}} &{\rm{in}}\quad  (\o-3h,-h)\\
\Du \quad  \tau_{-h}\ove{\what{\Phi}}\cp  
\tau_h W_o\cp\tau_{h}\ove{\what{\Phi}}\ &{\rm{in}} \quad(\o-2h,0)\\ 
\Du \quad \tau_{-2h}\ove{\what{\Phi}}\cp \tau_{-h} \ove{\what{\Phi}}\cp W_o &{\rm{in}} \quad (\o-h,h).\\
   \end{cases}
\end{equation}  
Since  $\tau_h W_o=W$ 
we obtain  \begin{equation}\label{maindual4:3}
 \what{\Phi^*}=\begin{cases} \Du\quad    \tau_{ h}W \cp\tau_{h}\ove{{\what{\Phi}}}\cp\tau_{2h} \ove{\what{\Phi}} &{\rm{in}}\quad  (\o-3h,-h)\\
\Du \quad  \tau_{-h}\ove{\what{\Phi}}\cp  
  W \cp\tau_{h}\ove{\what{\Phi}}\ &{\rm{in}} \quad(\o-2h,0)\\ 
\Du \quad \tau_{-2h}\ove{\what{\Phi}}\cp \tau_{-h} \ove{\what{\Phi}}\cp \tau_{-h}W  &{\rm{in}} \quad (\o-h,h).\\
   \end{cases}
\end{equation}  By comparing formulas (\ref{maindual4:2}) and  (\ref{maindual4:3}) we see that the expressions of $\what\Phi^*$ in  intervals symmetric with respect to zero can be obtained from each other   by replacing  $\tau_{jh}$ with $\tau_{-jh}$   (note that replacing $\tau_{jh}$ by $\tau_{-jh}$ changes also the sign of $W$). 
\par Formulas  (\ref{maindual4:1}), (\ref{maindual4:2}) and
  (\ref{maindual4:3})
give the  dual generators. 
 This completes the proof of the theorem when $\frac{\o}{2}< h< \frac{2\o}{3}$.
\par If $h=\frac{\o}{2}$  one argues as before; the only difference is that now the intervals $(\o-2h,2h-\o)$, $(-h,h-\o)$ and $(\o-h,h)$ are  empty.  \end{proof}
  Theorems   \ref{maindual3gen} and   \ref{maindual4gen} below generalize   Theorems \ref{maindual3} and  \ref{maindual4} respectively.  We omit the proofs, which are analogous to the proofs  of  Theorems \ref{maindual3} and  \ref{maindual4}. We recall that  $K_{-},K,$ and $K_{+}$ denote  the intervals $(0,\o-(\ell-1) h),$   $ (\o-(\ell-1) h,-\o+\ell h),$ and $(-\o+\ell h,h)$ 
  respectively. 
 \begin{theorem}\label{maindual3gen}
Let $\frac{\o}{\ell-\frac{1}{2}}  \le h <\frac{\o}{\ell-1}$ and denote by $\Phi$ the vector   $(\varphi_1,\varphi_2,\ldots,\varphi_{2\ell-1})$, where $\varphi_j\in \BL$, $j=1,\dots,2\ell-1.$ If  assumptions (\ref{maingendue:1})-(\ref{maingendue:3}) hold, i.e. if $\Efi$ is a frame for $\BL$, then  the Fourier transform of the   dual genera\-tors    $\Phi^*=(\varphi_1^*,\varphi_2^*,\ldots{\varphi}_{2\ell-1}^*)$ is 
  \begin{equation}\nonumber
\widehat{\Phi}^*= \begin{cases} 
(-1)^{\ell+1-k}  D  \quad \displaystyle{\crosspd{j=-\ell-k+1}{\ell-k-1}{j\not=0}\tau_{jh} \overline{\what\Phi}} 
  &{ \rm{in}}\quad \tau_{kh} K_-\quad
{\rm for}\   -(\ell-1)\le k\le\ell-1\\{\ }
  \\
  -\Du \ \displaystyle{ \crosspw{j=-\ell-k+1}{\ell-k-2}
  {\tau_{jh}  \overline{\what\Phi}}  { 
  \overline{\what\Phi}} {\tau_{-kh}W}}&{\rm{in}}
  \quad \tau_{kh} K\quad
{\rm for}\    -(\ell-1)
 \le k\le\ell-2
   \\ {\ }\\
 (-1)^{\ell+k} D  \qquad \displaystyle{\crosspd{j=-\ell-k}{\ell-2-k}{j\not=0}\tau_{jh} \overline{\what\Phi}}  &{\rm{in}} \quad
  \tau_{kh} K_+\quad
{\rm for}\    -\ell\le k\le\ell-2
  \end{cases}
  \end{equation}
where  $W=\displaystyle{\crossp{j=-\ell+1}{\ell-2}}\tau_{jh}\what{\Phi}$, $ \Du=\, h^{2\ell-3}\JJbstinv$   and  $D=h^{\ell-\frac{3}{2}}\,({\det}\   {\Jfistar})^{-1}$.  
 \end{theorem}
We recall that   
   $I_{-},I,I_{+}$ denote  the intervals $(0,-\o+\ell h),$ $ (-\o+\ell h,\o-(\ell-1)h),$ and   $(\o-(\ell-1)h,h)$  (see (\ref{I})). 
   \begin{theorem}\label{maindual4gen} 
 Let $\frac{\o}{\ell} \le h <\frac{\o}{\ell-\frac{1}{2}}$ and denote by $\Phi$   the vector  $(\varphi_1,\varphi_2, \ldots ,\varphi_{2\ell})$  where $\varphi_j\in \BL$, $j=1,\dots,2\ell.$ If  assumptions (\ref{maingen:1})-(\ref{maingen:3}) hold, i.e. if $\Efi$is a frame for $\BL$, then  the Fourier transform of the   dual generators    $\Phi^*=(\varphi_1^*,\varphi_2^*,\ldots{\varphi}_{2\ell}^*)$  is
 \begin{equation}\nonumber
\widehat{\Phi}^*= \begin{cases} 
 \ \Du \ \displaystyle{ \crosspw{j=-\ell-k+1}{\ell-k-1}{\tau_{jh} \ove{\what\Phi}}{\ove{\what\Phi}}{\tau_{-kh}W}  } 
  &{ \rm{in}}\quad \tau_{kh} I_- \quad
{\rm for}\quad   -(\ell-1)\le k\le\ell-1\\
\\
(-1)^{\ell+k} D  \hskip .1truecm \displaystyle{\crosspd{j=-\ell-k}{\ell-k-1}{j\not=0}\tau_{jh} \ove{\what\Phi}} 
  &{ \rm{in}}\quad  \tau_{kh}\, I   \qquad {\rm for}\   -\ell\le k\le\ell-1\\
    \\
  \Du \quad \displaystyle{\crosspw{j=-\ell-k}{\ell-k-2}{\tau_{jh} \ove{\what\Phi}}{\ove{\what\Phi}}{\tau_{-kh}W_o}  } &{\rm{in}} \quad\tau_{kh}\, I_{+}\quad {\rm for}\   -\ell\le k\le\ell-2 
  \end{cases}
  \end{equation}
  where     $W=\displaystyle{\crossp{j=-\ell+1}{\ell-1}}\tau_{jh}\what{\Phi},$ $\, W_o=\displaystyle{\crossp{j=-\ell}{\ell-2}}\tau_{jh}\what{\Phi}$, $\, \Du=h^{2\ell-2}\, \JJbstinv$ and  \break $D=h^{\ell-1}\,({{\det}}\ {\Jfistar})^{-1}.$ 
 \end{theorem}
 \section{Sampling formulas for the space $\BL $}
\label{s: Sampling}
   In this section   we shall  apply  the previous results    to  oversampling formulas for  the Hilbert transform sampling and the derivative sampling  in $\BL$.  
In the derivative sampling formula the  coefficients are the  values of the function and of its   derivatives $f^{(j)}$, $ 1\le  j \le  K,$  at  the points of  a uniform grid on $\R$.   It was  first  obtained    by D. Jagerman and L. Fogel for $K=1$ 
 and   by  Linden and N. M. Abramson for any $K$  \cite{JF}  \cite{L} \cite{LA}. Successively J. R. Higgins  derived  the same expansion formulas by using the Riesz basis method \cite{Hi}. \\ 
In \cite{SF} D.M.S. Santos and P.J.S.G. Ferreira    have   obtained  a two-channel derivative oversampling formula  for $B_{\omega_a}$ with $\omega_a<\omega$ by projecting both the   Riesz  basis generators  of the space $\BL$  and their  duals    into   the space  $B_{\omega_a}$. With this technique  the   projected family is a   frame; however   notice that     projecting  the  dual of a Riesz basis   does not yield the   dual frame. Thus the coefficients of the expansions of a function computed with respect to the projected duals   are not minimal in   least square norm.  
\par    Let $t_o$ be  such that $\o\le h<2\o$ and let  $\Phi=(\varphi_1,\varphi_2 )$ be a vector  such that $\Efi$ is a frame for $\BL.$ Then by (\ref{expansion})
\begin{equation}\label{Sformula}
f=\sum_{i=1,2}\sum_{k\in\Z} 
\langle f,\tau_{kt_o} \varphi_i\rangle \tau_{k t_o}\varphi_i^*
 \hskip 1truecm \forall f\in \BL.
\end{equation}
By using the  Plancherel and the inversion  formulas  we see that 
the coefficients
\begin{equation}\label{S:coeff}\langle f,\tau_{kt_o} \varphi_i \rangle=
\big({\cM}_j {f}\big)  (k t_o) \hskip 1truecm  j=1,2 
\end{equation}
are the samples of the functions    ${\cM}_j f =\cF^{-1}\what{\varphi}_j\,\cF f$ at the points  $kt_o, \, k\in \Z.$ For this reason (\ref{Sformula}) is called  a {\it sampling formula}. These formulas are useful in applications when one wants to reconstruct a signal from  samples taken from two transformed version of the signal. For instance one may want to reconstruct $f$ from samples of $f$ and $f'$ (derivative sampling)  or from samples of $f$ and  its Hilbert transform ${\mathcal H}f= -i\, {\mathcal F}^{-1}{\rm sign}{\mathcal F}f $ (Hilbert transform sampling).  Both are particular cases of the family of frames generated by the translates of two functions $\varphi_1,\varphi_2$ such that $\what{\varphi_1}=\chi_{[-\o,\o]}$, $\what{\varphi_2}=m  \chi_{[-\o,\o]},
$  where $m$ is a function in $L^\infty(\R)$.    
\begin{proposition}\label{m}     Let   $m$ be a function in $L^{\infty}(\R)$
 and let  $\Phi=(\varphi_1,\varphi_2)$ where 
 \begin{equation}\label{m:frame} 
 \what{\varphi_1}=\chi_{[-\o,\o]}  \hskip 1truecm   \what{\varphi_2}(x)=m  \chi_{[-\o,\o]}.
\end{equation}
Suppose that  $\o\le h< 2\o$.  Then   $\Efi$ is a frame for $\BL$  if and only if there exists a positive number $\eta$ such that 
\begin{equation}\label{detJm}
 |m-\tau_{-h}m|\ge \eta  \hskip1truecm   a.e.\, in  \ {(h-\o,\o)}.
\end{equation} The Fourier transforms of the dual generators are
 \begin{equation}\nonumber
\what{\varphi_1}^*=
\begin{cases}
 \frac{\tau_h \ove{m}}{h({\tau_h  \ove{m}}-\ove{m})} , & { \textnormal in }\,  [-\o,\o-h]\\
 {\, }\\
\frac{1}{h(1+|m|^2)},& { \textnormal in } \, (\o-h,h-\o) \\
{\, }\\
  \frac{-\tau_{-h}\ove{m}}{h({  \ove{m}}-\tau_{-h}\ove{m})} , &{ \textnormal in }  \, [h-\o,\o]\\ 
\end{cases}
\hskip .5truecm \what{\varphi_2}^*=
\begin{cases}
 \frac{-1}{h({\tau_h  \ove{m}}-\ove{m})} &  { \textnormal in }\,  [-\o,\o-h]\\
{\, }\\
\frac{m}{h(1+|m|^2)}, &  { \textnormal in }\,   (\o-h,h-\o) \\
{\, }\\
 \frac{1}{h({  \ove{m}}-\tau_{-h}\ove{m})} , &{ \textnormal in } \,   [h-\o,\o].\\
\end{cases}
\end{equation}
If $h=\o$ then $\Efi$ is a Riesz basis for $\BL$. 
\end{proposition}
\begin{proof} Since  $\det\, \Jfistar=h(\ove{m}-\tau_{-h}\ove{m})$ the assumptions of Theorem \ref{maingen} are satisfied.
The expression of  the Fourier transforms of  the dual generators can be easily obtained from  Theorem \ref{maindualnew}.
\end{proof} 
By choosing $m(x)=-i \, {\rm sign} (x)$    in (\ref{m:frame}) we obtain the Hilbert transform     frames  for $\BL$. For $h=\o$ the associated sampling formula is known as the {\it Hilbert transform sampling formula}   (see  \cite[Ex.12.9]{Hi}).  The coefficients of the expansion are the values of the function $f$  and   its Hilbert transform ${\mathcal H}f$ at the sample points $k t_o, k\in \Z$. Denote by $\snc $   the function $ ={\sin(x)}/{x}$. 
 \begin{corollary}\label{HTFRB}    Let 
     $ {\varphi_1}, {\varphi_2}$ be    defined by
  \begin{equation}\label{H:frame} 
 \what{\varphi_1}=\chi_{[-\o,\o]}  \hskip 1truecm   \what{\varphi_2}=-i  \chi_{[-\o,\o]} {\rm sign }.
\end{equation}
If $\o\le h<2\o$   then  $\Efi$ is a tigth frame for $\BL.$   The   dual generators are $ {\varphi^*_i}=
(2h)^{-1}  {\varphi_i}$ for $i=1,2.$ If $h=\o$ then $\Efi$  is a Riesz basis for $\BL.$
Moreover for any $f\in \BL$  the following Hilbert transform  sampling formula holds
\begin{equation}\nonumber 
f(x)=\frac{\o}{h} \sum_{k\in\Z} \Big(
f(k t_o) \tau_{-k t_o}\cos\big(\frac{\o x}{2}\big) \, \snc(\frac{\o x}{2})\ - \ (\mathcal{H} f)(k t_o) \tau_{-k t_o}\sin\big(\frac{\o x}{2}\big) \, \snc\big(\frac{\o x}{2}\big)  \Big).
\end{equation} 
\end{corollary}\begin{proof}
The assumptions of Proposition \ref{m} are   satisfied.   Thus $\Efi$ is a frame and the expression of the duals  follows immediately.
The     frame is tigth because $T T^*=  \frac{1}{2h}I$, since $\Phi^*=\frac{1}{2h}\Phi$. \par\noindent 
Standard calculations show that  if $\varphi_1$ and $\varphi_2$  are the functions  given by (\ref{H:frame}) then
 \begin{equation}\label{HTduals}\varphi_1(x)=\frac{\sqrt{2}}{\sqrt{\pi}}\o \cos\big(\frac{\o x}{2}\big) \, \snc\big(\frac{\o x}{2}\big)\hskip 1 truecm \varphi_2(x)=\frac{\sqrt{2}}{\sqrt{\pi}}\o \sin\big(\frac{\o x}{2}\big) \, \snc\big(\frac{\o x}{2}\big).  \end{equation}  
Moreover the    coefficients of the expansion   formula (\ref{Sformula})  are
 $  \sqrt{2\pi}\hat{f}(k t_o) $ and \break $ -\sqrt{2\pi}(\mathcal{H}f)(kt_o)$, $k\in \Z$. This  proves  also the expansion formula.
\end{proof}
    \par\noindent
   By choosing    $m(x)=ix$   in (\ref{H:frame})  we obtain the derivative  frame  for $\BL$.  Given a function $g$ we shall denote by $g_{\delta}$ the  function $  g(\delta x).$
 Note that $\delta\, \snc_{\delta}= \sqrt{\frac{\pi}{2}}  \what{ \chi}_{_{[-\delta,\delta]}}.$
 \begin{corollary}\label{DFRB}   Let $ {\varphi_1}, {\varphi_2}$ be defined by 
 \begin{equation}\label{Deriv:frame} 
 \what{\varphi_1}=\chi_{[-\o,\o]}  \hskip 1truecm   \what{\varphi_2}=i x   \chi_{[-\o,\o]} . 
\end{equation}
If $\o\le h<2\o$   then  $\Efi$ is a   frame for $\BL$;  if $h=\o$ then it is a Riesz basis for $\BL.$ The Fourier transforms of the dual generators are 
 \begin{equation} \nonumber 
\what{\varphi_1}^*(x)=
\begin{cases}
\frac{1}{h}(1-  \frac{|x|}{h})  &\hskip-.1truecm  h-\o<|x|<\o\\
   \frac{1}{h (1+x^2)} &\hskip-.1truecm|x|<h-\o  \\
\end{cases}
\hskip 0.4truecm \what{\varphi_2}^*(x)=
\begin{cases}
\frac{i}{h^2} \,{\rm sign}(x)  &\hskip-.1truecm h-\o<|x|<\o\\
  \frac{ix}{h (1+x^2)}  &\hskip-.1truecm |x|<h-\o.  \\
  \end{cases}
\end{equation} 
Moreover for  any $f\in \BL$ the following  derivative sampling  formula holds
 \begin{equation}\label{dsamp:expdue}
f=\sqrt{2\pi}\sum_{k\in \Z} \Big( f(kt_o)\tau_{-kt_o}\varphi_1^*- f'(kt_o) \tau_{-kt_o}\varphi_2^*\Big),
\end{equation}
 where    
 \begin{equation}\nonumber
 \begin{aligned}
\varphi_1^*(x)=&e^{-|\cdot|} \ast(h-\o)\frac{1}{h}\snc_{h-\o}  (x)
+\frac{1}{\pi}\big(\o\, \snc_{\o}-(h-\o)\,\snc_{h-\o}\big)\ast \snc^2_{h/2}(x)  \\
{\ }\\
\varphi_2^*(x)=&e^{-|\cdot|}\ast  (h-\o)\frac{1}{h}\snc^{'}_{h-\o} (x)+
\frac{\sqrt{2}}{\sqrt{\pi}}\frac{1}{ x h^2}  \big( \cos_{\o}-\cos_{h-\o}\big)(x).
\end{aligned}\end{equation}
\end{corollary}\begin{proof}
By Proposition \ref{m}  $\Efi$ is a frame for all $\o\le h< 2 \o$ and it is a Riesz basis if $h=\o$.   From  (\ref{S:coeff}) we see  that  the coefficients of the expansion  (\ref{Sformula})  are $ \sqrt{2\pi}\hat{f}(k t_o)$  and  $-\sqrt{2\pi}f'(kt_o).$    The  expression   of the dual generators  
   can be obtained  from $\what{\varphi_1}^*$ and $\what{\varphi_2}^*$ by computing the inverse Fourier transform.
   \end{proof}
A simple calculation shows that if $h=\o$ then
$$
\varphi_1^*(x)=\frac{1}{\sqrt{2\pi}}\  {\snc}^2(\frac{\o x}{2}) \qquad \varphi_2^*(x)=- \frac{1}{ \sqrt{2 \pi}} \ x \,{\snc}^2(\frac{\o x}{2}).
$$
 Figure 1  and Figure 2 below show the  Fourier transforms of  the dual generators   in the Riesz basis case $h=\o=1$,   and in the case 
 $\o=1$ and $h=\frac{3}{2}\o$  respectively.
\begin{figure}[htbp]
\begin{center}
\begin{minipage}[c]{.40\textwidth}
\centering\setlength{\captionmargin}{0pt}
\includegraphics[width=5cm,height=6cm] 
 {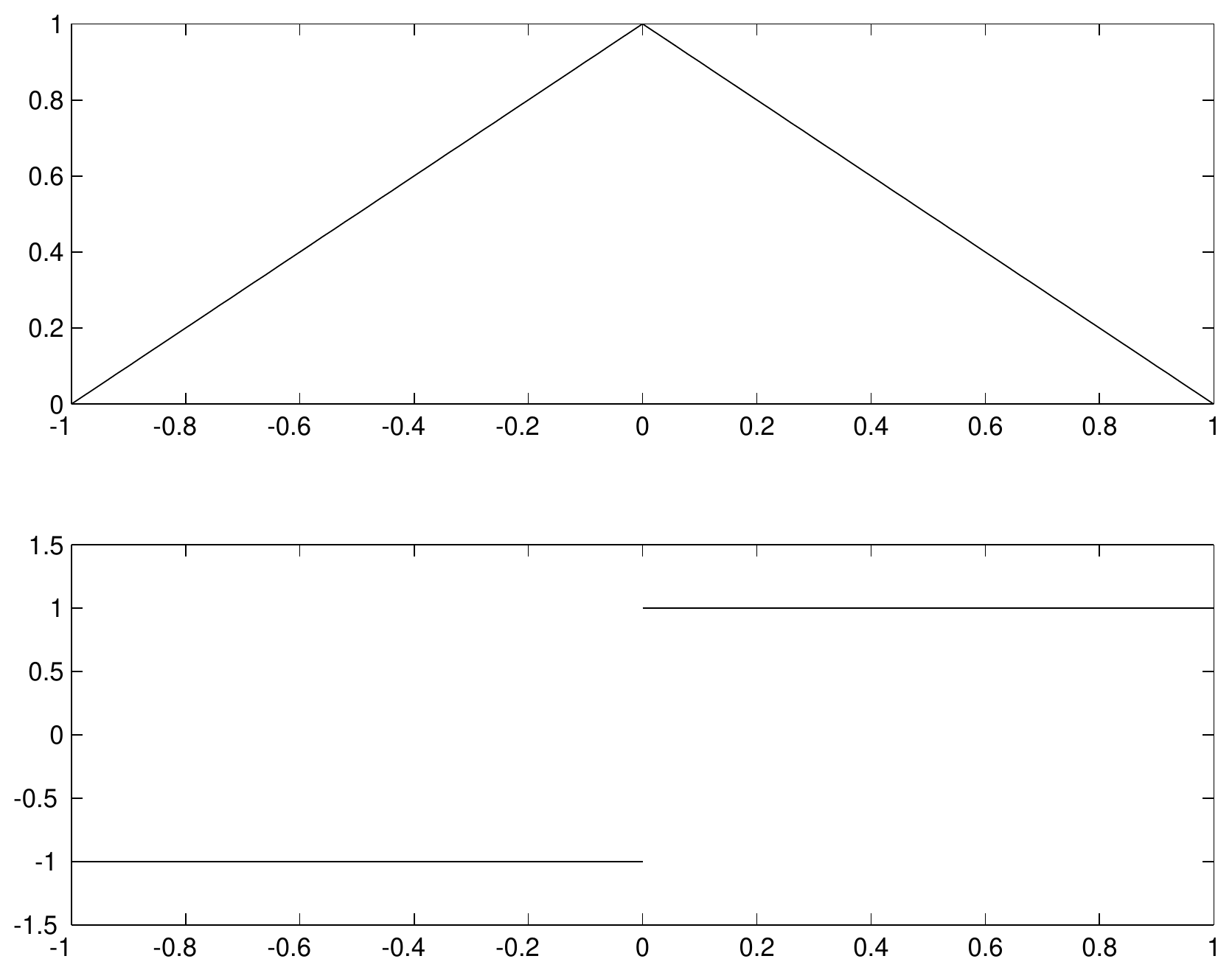}
\caption{$ \widehat{\varphi_1^{\star}},\frac{1}{i}\widehat{\varphi_2^{\star}} $ for $h=\o=1$}
\end{minipage}
\begin{minipage}[c]{.40\textwidth}
\centering\setlength{\captionmargin}{0pt} 
\includegraphics[width=5cm,height=6cm]
{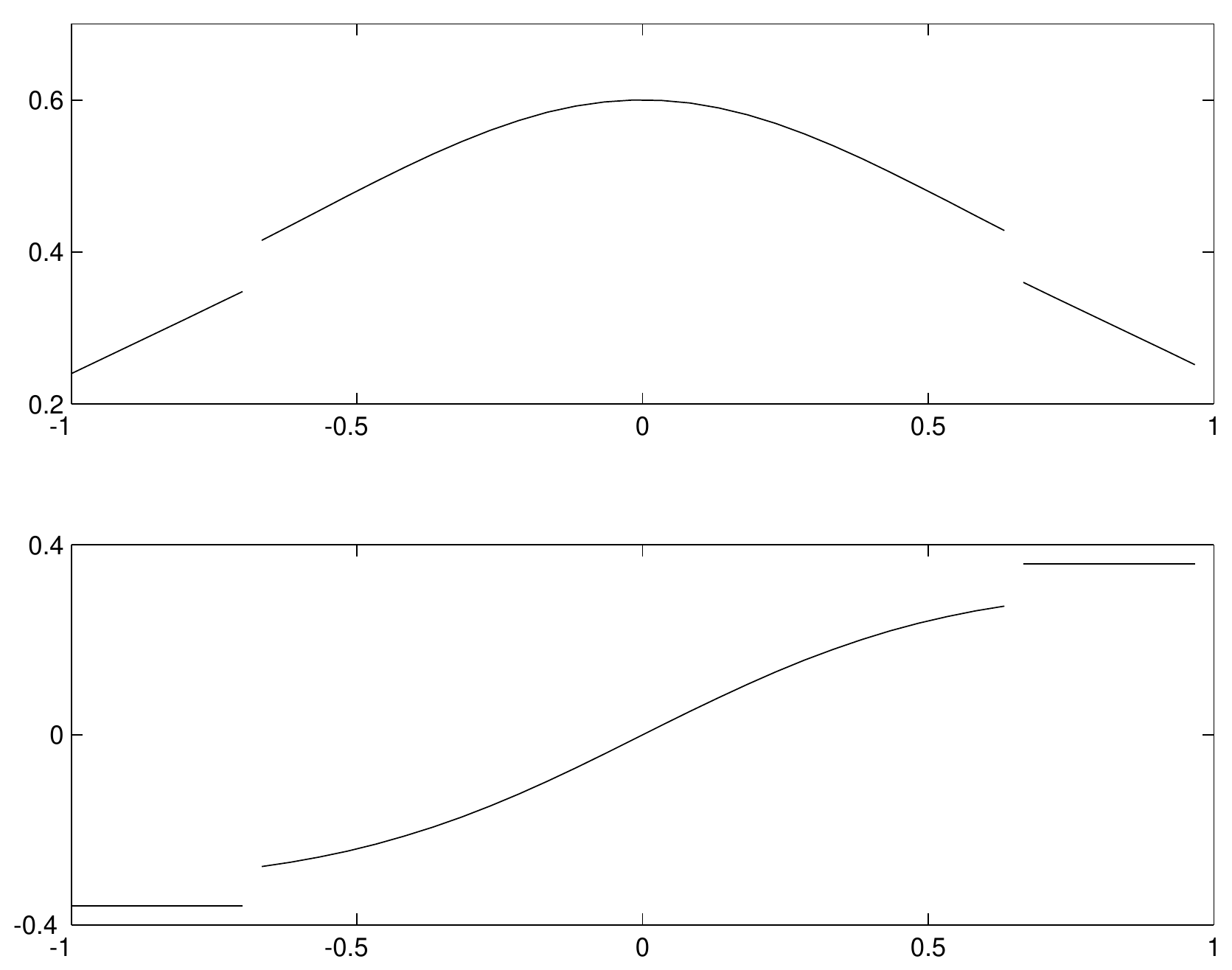}
\caption{$ \widehat{\varphi_1^{\star}},\frac{1}{i}\widehat{\varphi_2^{\star}} $ for $h=\frac{3}{2}\o=\frac{3}{2}$}
\end{minipage}
\end{center}
\end{figure}
 \par
{\bf Remark}  In   \cite{JF}  Jagerman and Fogel proved the following  two-channel derivative  sampling formula  in the case $h=\o$ i.e. $t_o=\frac{2\pi}{\o}$ 
 \begin{equation}\nonumber
f(x)=\sum_kf(kt_o) \Big(\snc^2\frac{\o}{2}\big(x-kt_o\big)\, +\, \frac{2}{\o} f'(kt_o)\sin \frac{\o}{2}\big(x-kt_o\big)\ \snc\frac{\o}{2}\big(x-kt_o\big)\, \Big)
\end{equation}
   (see   also \cite[p.135]{Hi}). Thus formula  (\ref{dsamp:expdue}) is  an extension of the case $t_o=\frac{2\pi}{\o}$   to   all   values of $t_o\in [\frac{\pi}{\o},\frac{2\pi}{\o})$ (i.e. for all $h$ 
such that $\o\le h< 2\o$). 

Our last example  is a  three channel derivative oversampling    formula. To obtain a frame  with three generators we must choose $\frac{2}{3}\o\le h<\o$ i.e. $t_o\in  \Big[\frac{3\pi}{\o},\frac{2\pi}{\o}\Big)$.  \par \noindent 
   \begin{corollary}\label{DFRB3}   Let $\what{\varphi_1},\what{\varphi_2},\what{\varphi_3}$ be defined by 
 \begin{equation}\nonumber
 \what{\varphi_1}=\chi_{[-\o,\o]},  \hskip 1truecm   \what{\varphi_2}=i x   \chi_{[-\o,\o]},  \hskip 1truecm   \what{\varphi_3}= - x^2   \chi_{[-\o,\o]}. 
\end{equation}
If $\frac{2\o}{3}\le h<\o$   then  $\Efi$ is a   frame for $\BL$;  if $h=\frac{2\o}{3}$ then it is a Riesz basis for $\BL.$ The Fourier transform of the dual generators are 
\begin{equation} \label{Der:genduali3} 
 \what{\Phi^{\ast}}(x)=
\begin{cases} 
\frac{1}{2h^3}\big( x^2+3hx+2h^2 ,-i (2x+3h),-1\big) \quad&  -\o<x<\o-2h\\
\big(A_h(x),B_h(x),C_h(x)\big) &\o-2h<x<h-\o \\
 \frac{1}{h^3} \big(h^2-x^2,2ix,1\big)   & h- \o<x<\o-h\\
 \big(A_{-h}(x),B_{-h}(x),C_{-h}(x)\big)   &\o-h<x<2h-\o\\
  \frac{1}{2h^3}\big( x^2-3hx+2h^2,-i (2x-3h),-1\big)  & 2h-\o<x<\o\\
\end{cases}
 \end{equation} 
where $$A_h(x)=  \frac{1}{h^2}\  \frac{(x+h) +(2x+h)(x+h)^2} {1+(2x+h)^2+x^2(x+h)^2}\hskip .6truecm B_h(x)=  \frac{-i}{h^2}\  \frac{1-x(x+h)^3} {1+(2x+h)^2+x^2(x+h)^2 }$$$$ C_h(x)=  \frac{1}{h^2}\  \frac{h+2x+x(x+h)^2} {1+(2x+h)^2+x^2(x+h)^2} .$$\\
Note that if  $h=\frac{2}{3}\o$    the intervals $(\o-h,2h-\o)$ and $(2h-\o,h-\o)$ are empty.  
\end{corollary}
 \begin{proof}
The family $\Efi$ is a frame for $\BL$ by      Theorem \ref{maindual3}. To compute  the dual generators we used   formula  (\ref{maindual3t}). 
 \end{proof} 
Thus    for each value of the parameter $h$ in $ [\frac{2 }{3}\o,\o)$ the family $\Efi$ is a frame and every signal in $\BL$ can be reconstructed from the values $f(kt_o),f^{(1)}(kt_o),f^{(2)}(kt_o), k\in\Z,$ by   the following three-channel derivative sampling formula  
 \begin{equation}\label{DSformula}
f=\frac{1}{\sqrt{2\pi}}\sum_{i=1}^{3}\sum_{k\in\Z} 
(-1)^{i-1} f^{(i-1)}(kt_o) \varphi_i^*(x-kt_o),
\end{equation}
where the dual  generators    are  given by (\ref{Der:genduali3}).\par \noindent
In the rest of this section we present two numerical experiments of reconstruction of a band limited signal by using  formula (\ref{DSformula}). We have chosen the   signal  $f(x)=\frac{1}{\sqrt{2\pi}} (\snc(x/2))^2$ (see Figure 3); note that the function $f$  is the Fourier transform of the function  $(1-|x|)_+$, therefore $\o=1$. \begin{figure}[h]
\begin{center}
\begin{minipage}[c]{.90\textwidth}
\centering\setlength{\captionmargin}{0pt}%
\includegraphics[width=6.5cm,height=4cm]{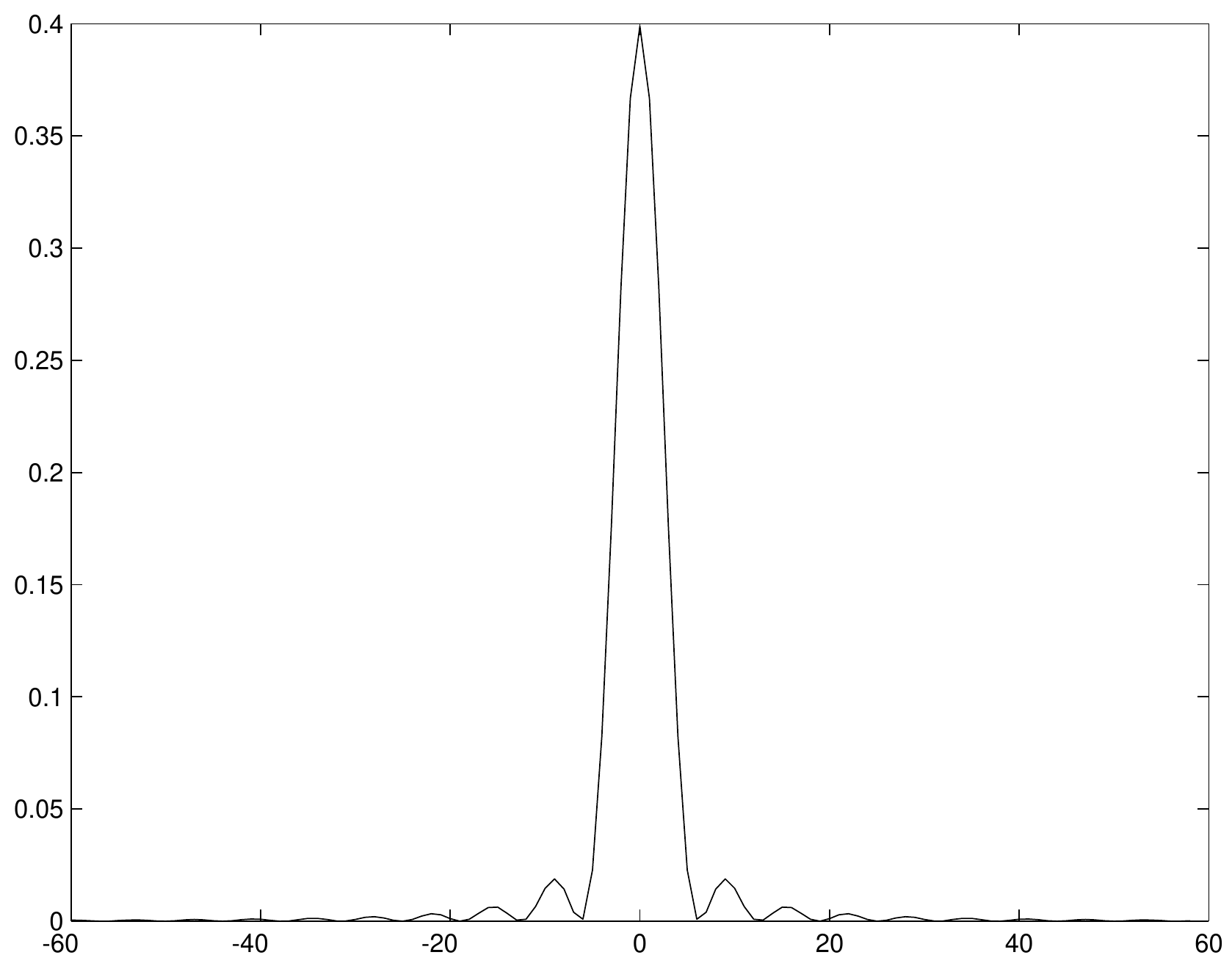}
\caption{Signal $f(x)=\frac{1}{\sqrt{2\pi}}(\snc(x/2))^2$}
\end{minipage}%
 \end{center}
\end{figure}
\par
In the first experiment we take    $h=\frac{2}{3}\o$   so that $\Efi$ is  a Riesz basis; in the second experiment  we take $h=\frac{11}{15}\o$ so that $\Efi$ is a frame. We observe that  in  the first case    it is possible to find  the analitic expression  of the functions $ \varphi_j^*$    in (\ref{Der:genduali3}), while in the second case they must be computed numerically. Indeed for  $h=\frac{2\o}{3}$   one obtains
 \begin{align}\nonumber
 \varphi_1^*&=\frac{1}{\sqrt{2\pi}}\Big(1+\frac{\o^2 x^2}{18}\Big)\snc^3\Big(\frac{\o}{3}x\Big) \nonumber
\\ \nonumber
\varphi_2^*&=-\frac{1}{\sqrt{2\pi}}x\, \snc^3\Big(\frac{\o}{3}x\Big)
 \\ \varphi_3^*&=\frac{1}{2\sqrt{2\pi}}x^2\, \snc^3\Big(\frac{\o}{3}x\Big).\nonumber
\end{align} 
 Figures 4 and  5 below show the functions $  \varphi_1^*, \varphi_2^*, \varphi_3^*$ and $ \widehat{ \varphi_1^*},\frac{1}{i}  \widehat{\varphi_2^*},  \widehat{\varphi_3^*}$ for $h=\frac{2}{3}$. 
 \begin{figure}[h]
\begin{center}
\begin{minipage}[c]{.40\textwidth}
\centering\setlength{\captionmargin}{0pt}%
\includegraphics[width=5cm,height=6cm]
{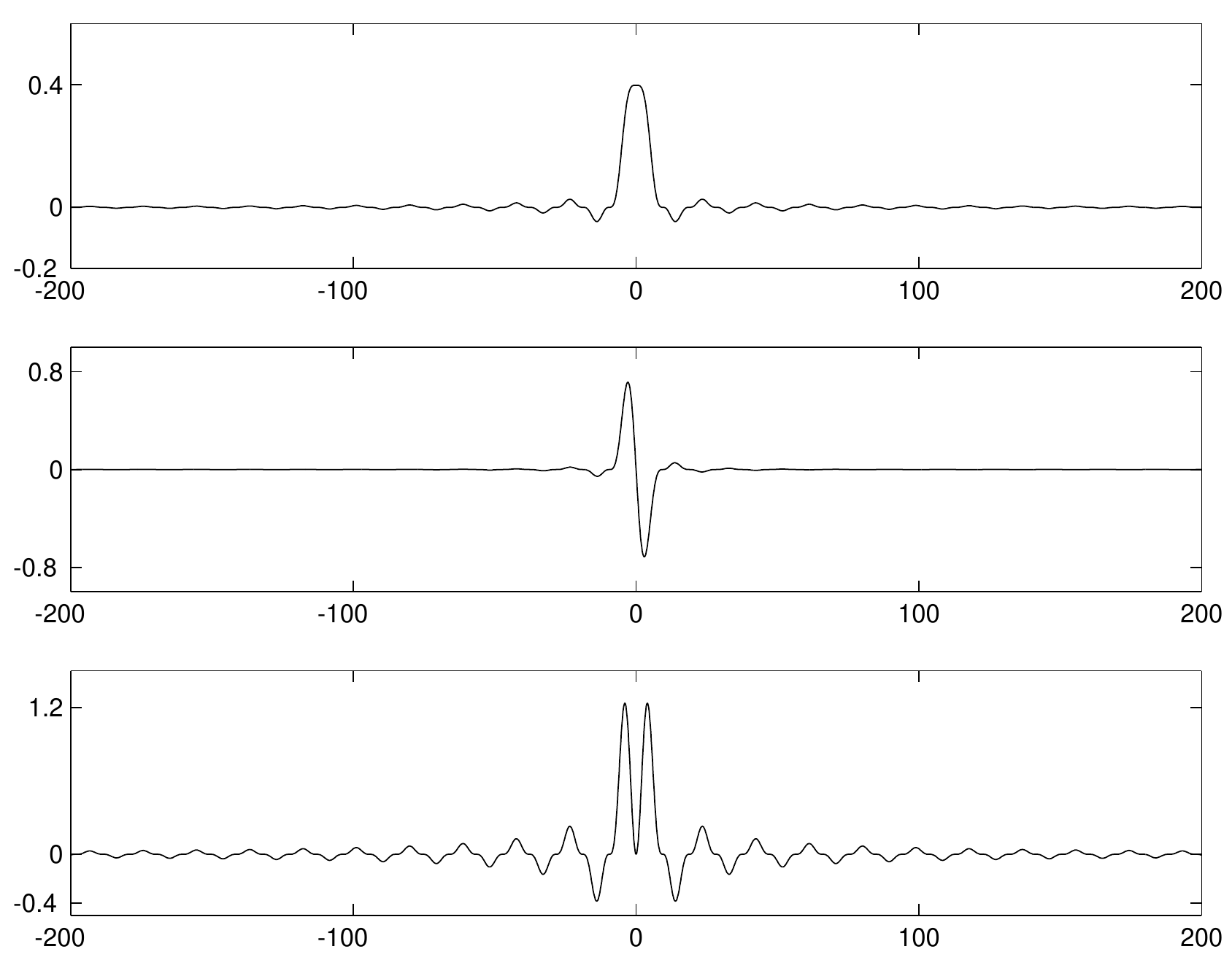}
\caption{$ {\Phi}^{\star};$  $h=\frac{2}{3}\o$}
\end{minipage}
\begin{minipage}[c]{.40\textwidth}
\centering\setlength{\captionmargin}{0pt}%
\includegraphics[width=5cm,height=6cm]
{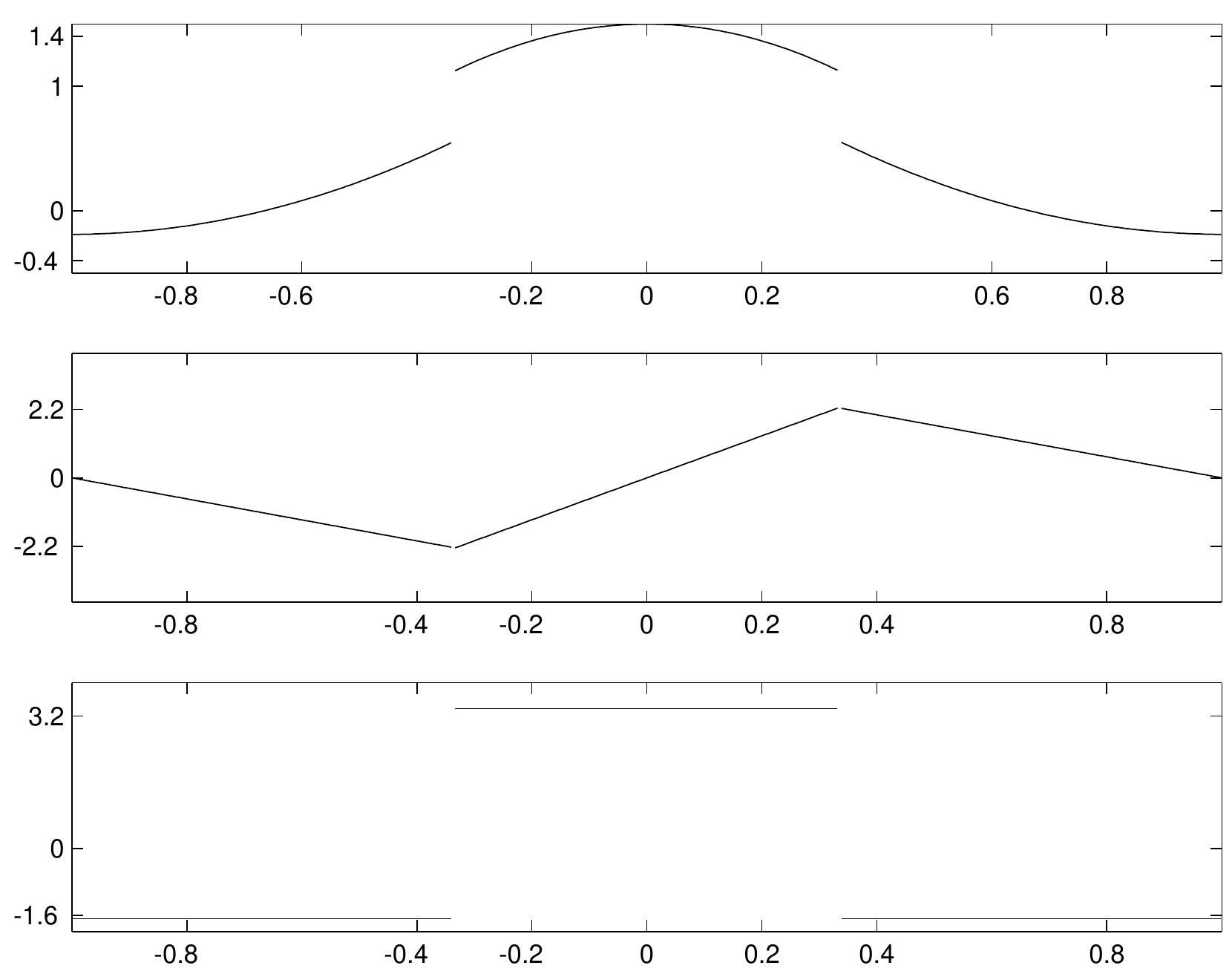}
\caption{$ \widehat{ \varphi_1^*},\frac{1}{i}  \widehat{\varphi_2^*},  \widehat{\varphi_3^*};$ $h=\frac{2}{3}\o$}
\end{minipage}
\end{center}
\end{figure}
 \par
Since  in the case  $h=\frac{2}{3}\o$ the expression of the duals is known,   it is possible to write explicitly the 
 sampling formula  (\ref{DSformula}):
  \par
\begin{align} \nonumber f(x)= &\frac{1}{\sqrt{2\pi}} \sin^3\big(\frac{\o x}{3}\big) \sum_n\,  \big(-1\big)^n \Big[\,
  f\big(\frac{3n\pi}{\o}\big) \frac{1}{\big(x- \frac{3n\pi}{\o}\big)^3} +
  \\ &\nonumber
  \\ &
  \nonumber
  +  f\big(\frac{3n\pi}{\o}\big)  \frac{ {\o^2}/{9}}{\big(x- \frac{3n\pi}{\o}\big)}- f^{(1)}\big(\frac{3n\pi}{\o}\big) \frac{1}{\big(x- \frac{3n\pi}{\o}\big)^2}+   f^{(2)}\big(\frac{3n\pi}{\o}\big)  \frac{1}{2\big(x- \frac{3n\pi}{\o}\big)^2}\,\Big].
\end{align}
 This formula was   first   given by Linden in \cite{L} and Linden and Abramson in \cite{LA}, see also \cite{HS}.   \par
In the case $h=\frac{11}{12}\o$  the dual generators   have been  obtained by computing numerically the inverse Fourier transforms of the functions in  (\ref{Der:genduali3}):
Figure 6 and Figure 7   below show the functions $  \varphi_1^*, \varphi_2^*, \varphi_3^*$ and $ \widehat{ \varphi_1^*},\frac{1}{i}  \widehat{\varphi_2^*},  \widehat{\varphi_3^*};$   $h=\frac{11}{15}\o.$     \begin{figure}[h]
\begin{center}
\begin{minipage}[c]{.40\textwidth} 
\centering\setlength{\captionmargin}{0 pt}%
\includegraphics[width=5cm,height=6cm]{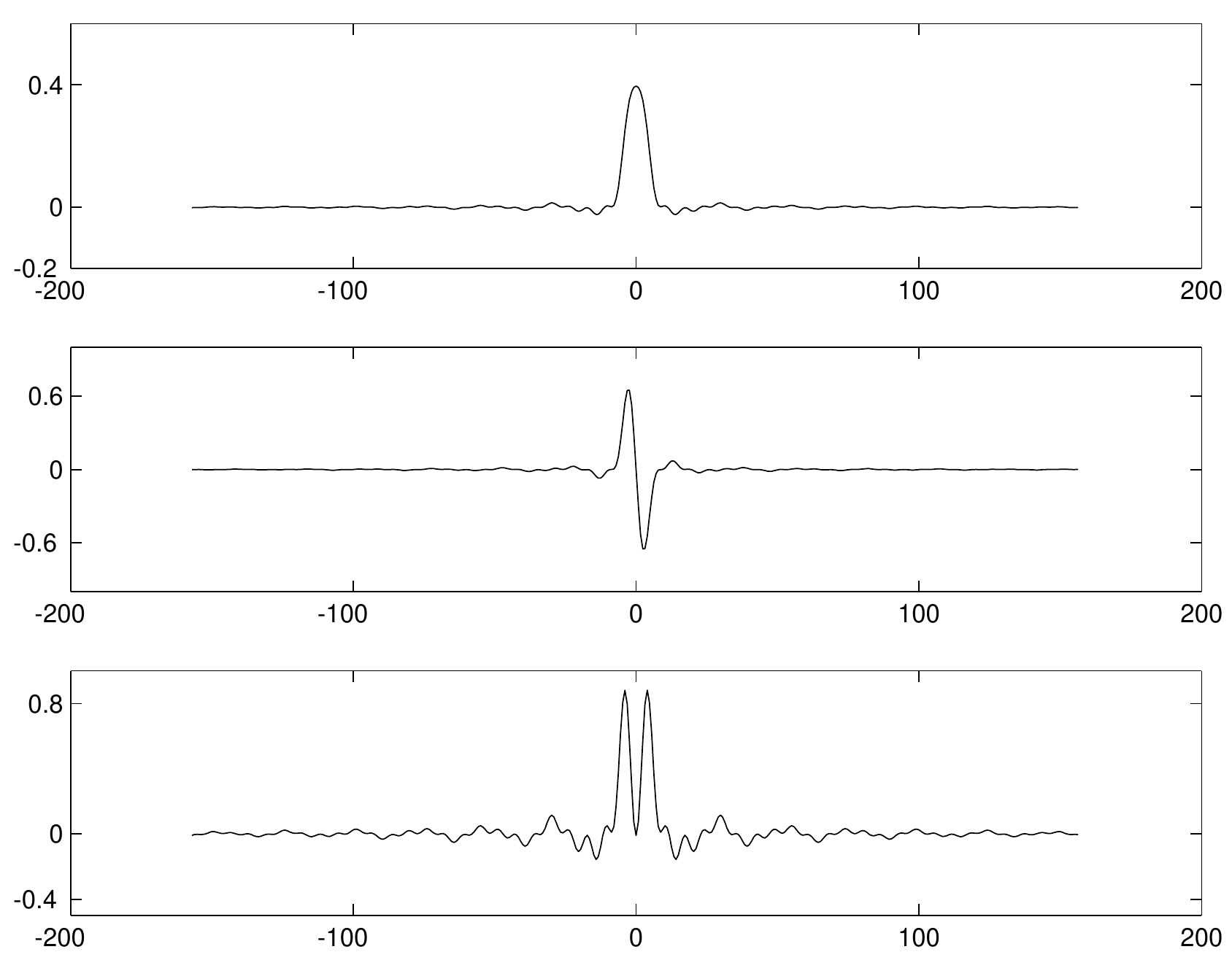}
\caption{${\Phi}^{\star}$ for  $h=\frac{11}{15}\o$}
\end{minipage}
\begin{minipage}[c]{.40\textwidth}
\centering\setlength{\captionmargin}{0pt}%
\includegraphics[width=5cm,height=6cm]{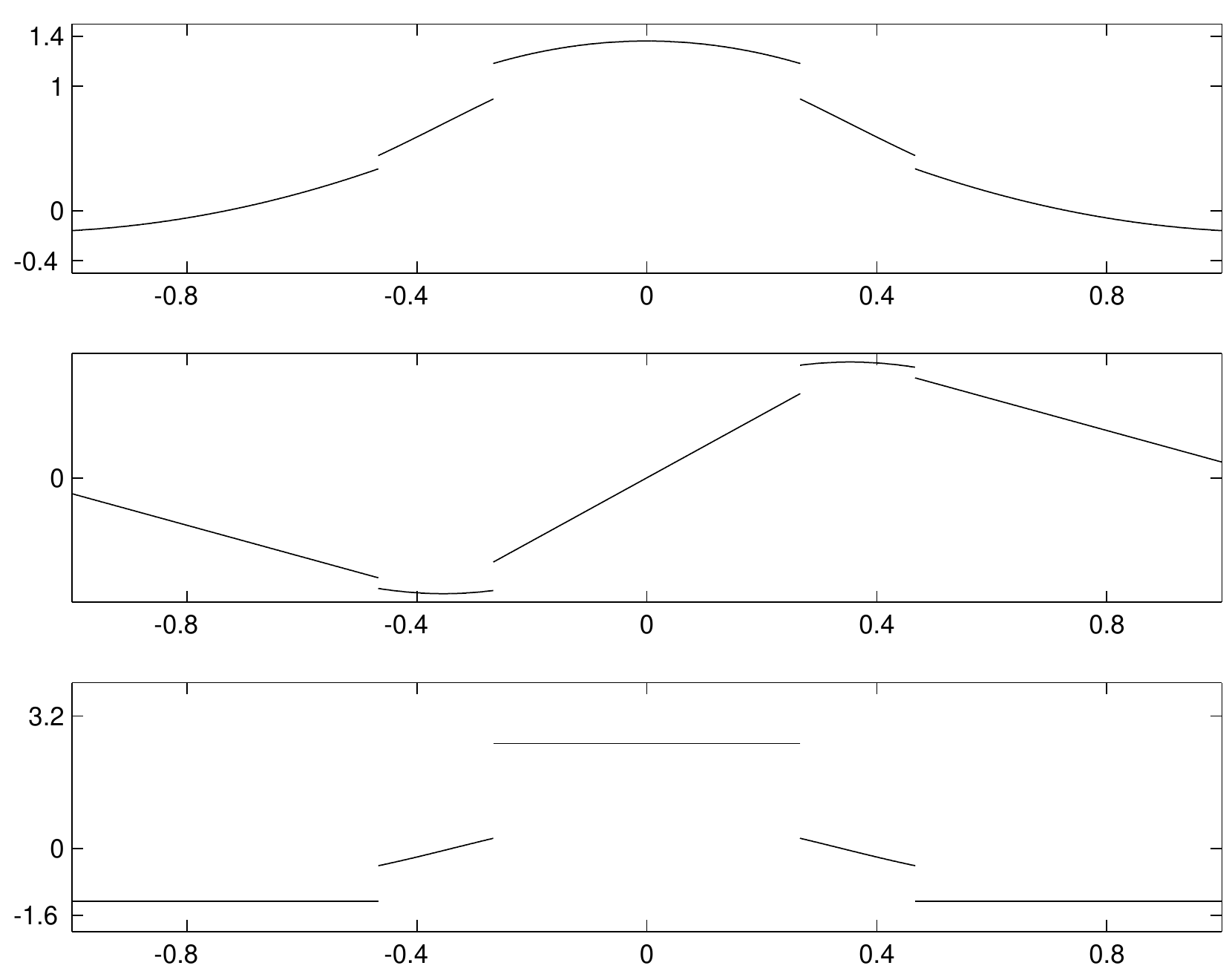}
\caption{$ \widehat{ \varphi_1^*},\frac{1}{i}  \widehat{\varphi_2^*},  \widehat{\varphi_3^*};$ $h=\frac{11}{15}\o$}
\end{minipage}
\end{center}
\end{figure}\par
  Figure 8 and  Figure 9 show  the error in  the  cases $ h=\frac{2}{3}\o$ and   $h=\frac{11}{15}\o$, respectively. Notice that in the first case the order of magnitude of the error is $10^{-4}$ while in the second case is $10^{-3}$.
In  second case, since the functions $\varphi_j^*$ were computed  numerically, to compute their values at the points $x-kt_o$ we used spline interpolation. This accounts for the different order of magnitude of the error.
\begin{figure}[h]
\begin{center}
\begin{minipage}[c]{.40\textwidth}
\centering\setlength{\captionmargin}{0pt}%
\includegraphics[width=.80\textwidth]{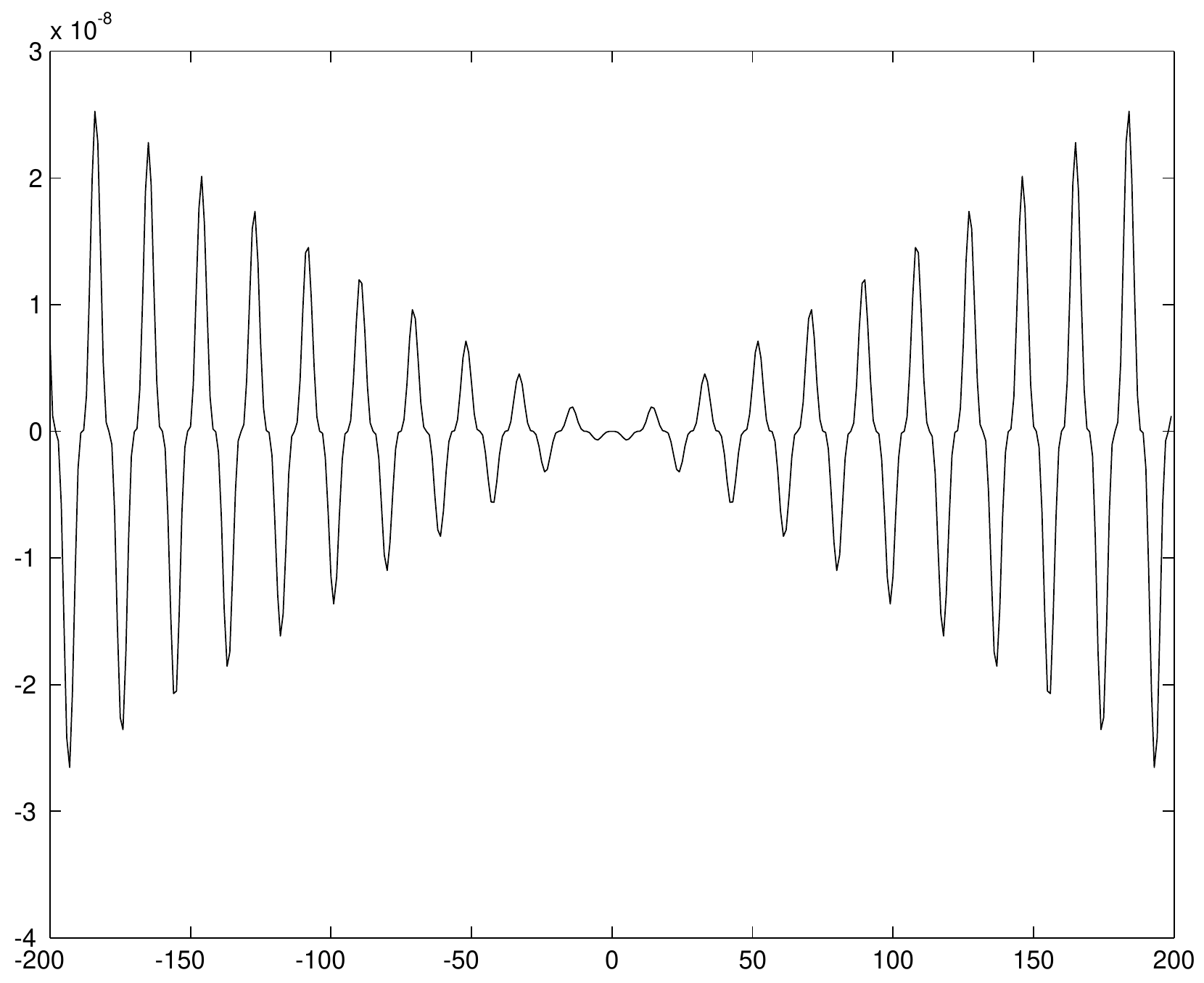}
\caption{The error  with
 $h=\frac{2}{3}$ }
\end{minipage}
\begin{minipage}[c]{.40\textwidth}
\centering\setlength{\captionmargin}{0pt}%
\includegraphics[width=.80\textwidth]{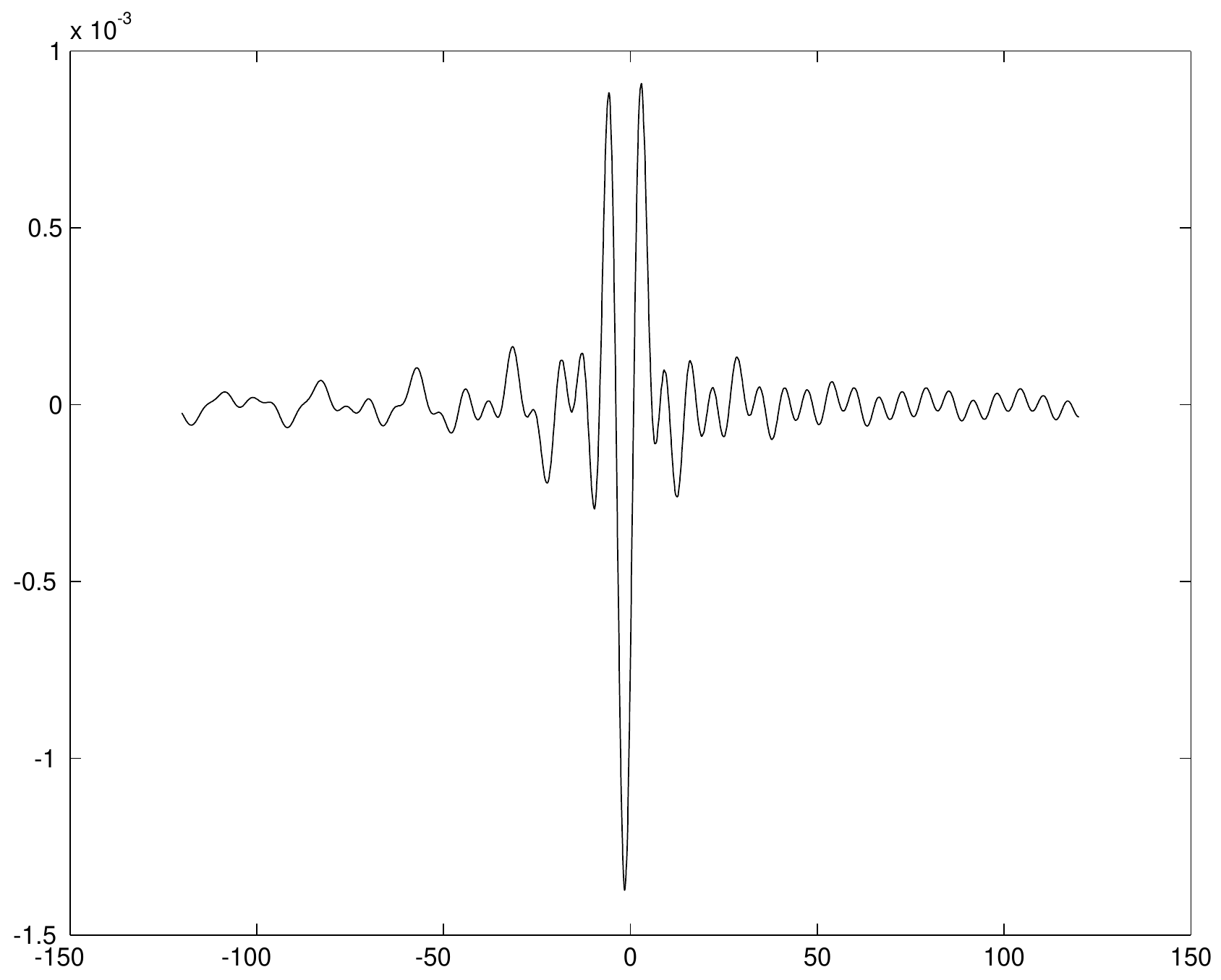}
\caption{The error with 
$h=\frac{11}{15}$ }
\end{minipage}
\end{center}
\end{figure}
 \vskip.5truecm
\centerline{ACKNOWLEDGMENTS}
\vskip.1truecm
The author thanks Aldo Conca  for   helpful discussions on the Moore-Penrose inverse of a matrix.

   \end{document}